\documentclass[11pt,english]{article}
\usepackage{lmodern}
\usepackage{lmodern}
\usepackage[T1]{fontenc}
\usepackage[latin9]{inputenc}
\usepackage{geometry}
\usepackage{color}
\usepackage{amsmath}
\usepackage{amssymb}

\makeatletter

\providecommand{\tabularnewline}{\\}

\newenvironment{lyxlist}[1]
{\begin{list}{}
{\settowidth{\labelwidth}{#1}
 \setlength{\leftmargin}{\labelwidth}
 \addtolength{\leftmargin}{\labelsep}
 }}
{\end{list}}

\usepackage{babel}

\date{June 28, 2019}

\usepackage{babel}

\makeatother

\usepackage{babel}
\begin{document}

\title{Statistical inference and feasibility determination: a nonasymptotic
approach\thanks{Two earlier versions (arXiv:1808.07127v1 and arXiv:1808.07127v2) of
this paper were circulated under the titles ``Concentration based
inference in high dimensional generalized regression models (I: statistical
guarantees)'' and ``Concentration based inference for high dimensional
(generalized) regression models: new phenomena in hypothesis testing''.
Substantial changes have been made in this version (arXiv:1808.07127v3).
This version was prepared during my appointment at Department of Statistics
and Department of Computer Science, Purdue University, West Lafayette,
IN. A start-up fund from Purdue College of Science supported this
research. Earlier versions were prepared during my appointment at
Michigan State University (Department of Economics, Social Science
Data Analytics Initiative) that also provided financial support. }}

\author{Ying Zhu\thanks{Email: yiz012@ucsd.edu. Assistant Professor of Economics (starting
on July 1, 2019), Department of Economics, University of California,
San Diego. }}
\maketitle
\begin{abstract}
We develop non-asymptotically justified methods for hypothesis testing
about the $p-$dimensional coefficients $\theta^{*}$ in (possibly
nonlinear) regression models. Given a function $h:\,\mathbb{R}^{p}\mapsto\mathbb{R}^{m}$,
we consider the null hypothesis $H_{0}:\,h(\theta^{*})\in\Omega$
against the alternative hypothesis $H_{1}:\,h(\theta^{*})\notin\Omega$,
where $\Omega$ is a nonempty closed subset of $\mathbb{R}^{m}$ and
$h$ can be nonlinear in $\theta^{*}$. Our (nonasymptotic) control
on the Type I and Type II errors holds for fixed $n$ and does not
rely on well-behaved estimation error or prediction error; in particular,
when the number of restrictions in $H_{0}$ is large relative to $p-n$,
we show it is possible to bypass the sparsity assumption on $\theta^{*}$
(for both Type I and Type II error control), regularization on the
estimates of $\theta^{*}$, and other inherent challenges in an inverse
problem. We also demonstrate an interesting link between our framework
and Farkas' lemma (in math programming) under uncertainty, which points
to some potential applications of our method outside traditional hypothesis
testing.
\end{abstract}

\section{Introduction}

A common feature of the existing procedures that are deemed ``practical''
for testing nonlinear hypotheses about regression coefficients is
that they hinge on asymptotic validity to some extent. This occurrence
is perhaps not coincidental as asymptotic analyses often allow one
to focus on some ``leading'' term(s) by assuming the ``remainder''
term(s) approach to zero faster, which can be quite convenient for
determining the threshold in a test. However, many real-world applications
(in controlled experiments, for example) have a limited sample size
which renders any asymptotic argument questionable. Our primary goal
is to find situations where effective non-asymptotic methods can be
developed for hypothesis testing about the coefficients $\theta^{*}\in\mathbb{R}^{p}$
in regression models. Given a function $h:\,\mathbb{R}^{p}\mapsto\mathbb{R}^{m}$,
let
\begin{equation}
H_{0}:\,h(\theta^{*})\in\Omega\;\textrm{vs.}\;H_{1}:\,h(\theta^{*})\notin\Omega,\label{eq:nonlinear}
\end{equation}
where $\Omega$ is a nonempty closed subset of $\mathbb{R}^{m}$ and
$h$ is allowed to be nonlinear in $\theta^{*}$. Relative to existing
literature, we consider these broader forms of hypotheses and the
impact of $m$, the number of restrictions in the null hypothesis
$H_{0}$. Throughout this paper, we assume that $\left\{ \theta\in\mathbb{R}^{p}\,:\,h(\theta)\in\Omega\right\} \neq\emptyset$
and $H_{0}$ does not contain any redundant restrictions. 

Our main focus is the following Gaussian regression model
\begin{equation}
Y_{i}=g\left(V_{i};\theta^{*}\right)+W_{i},\qquad i=1,...,n.\label{eq:reg}
\end{equation}
In the equation above, the functional form of $g\left(V_{i};\theta^{*}\right)$
is known and may be nonlinear in $\theta^{*}$; $Y=\left\{ Y_{i}\right\} _{i=1}^{n}$
is an $n-$dimensional vector of responses; $V=\left\{ V_{i}\right\} _{i=1}^{n}\in\mathbb{R}^{n\times k}$
is the matrix of covariates with the $i$th row specified by $V_{i}$;
$W=\left\{ W_{i}\right\} _{i=1}^{n}\sim\mathcal{N}(\mathbf{0}_{n},\,\sigma^{2}\boldsymbol{I}_{n})$
and is independent of $V$, where $\mathbf{0}_{n}$ denotes an $n-$dimensional
vector of zeros; $\theta^{*}$ is a $p-$dimensional vector of unknown
coefficients and $p$ is allowed to exceed $k$ as well as the sample
size $n$. Throughout the paper, we make our argument by conditioning
on $V$.

In (\ref{eq:nonlinear}), allowing $h$ to be nonlinear functions
of $\theta^{*}$ can be very useful. For example, policy researchers
are often interested in testing whether $\frac{1}{n}\sum_{i=1}^{n}\frac{\partial g\left(v_{i};\theta^{*}\right)}{\partial v_{ij}}$
or $\frac{1}{n}\sum_{i=1}^{n}\frac{\Delta g\left(v_{i};\theta^{*}\right)}{\Delta v_{ij}}$
lies in some interval. These quantities are referred to as the average
partial effect (APE) of $v_{ij}$s on $\mathbb{E}\left(Y_{i}\vert V_{i}=v_{i}\right)=g\left(v_{i};\theta^{*}\right)$,
holding all the other $v_{il}$s constant ($l\neq j$). For evaluating
policy interventions in social science or treatment procedures in
medical studies, both applied and theoretical researchers have recognized
the importance of APEs for understanding the magnitudes of effects
(see \cite{w10} for more discussions on this topic). If $g\left(v_{i};\theta^{*}\right)=v_{i}\theta^{*}$
in (\ref{eq:reg}), the APE of $v_{ij}$s is simply $\theta_{j}^{*}$.
However, when a nonlinear model is used to analyze the effects of
$v_{ij}$s, the APEs often depend on $\theta^{*}$ in a nonlinear
fashion and the individual coefficients themselves may no longer convey
interpretable information about the effects of $v_{ij}$s. 

This paper makes several contributions. First, the new method we propose
provides a finite sample alternative to the classical asymptotic procedures
(such as the Wald, Lagrange multiplier, and likelihood ratio tests)
for testing a single or multiple nonlinear hypotheses about coefficients
in regression models where $p$ is small relative to $n$. Second,
this method can be used for testing simultaneous (nonlinear) hypotheses
when $p$ is comparable to or larger than $n$ and the number of restrictions
in $H_{0}$ is large relative to $p-n$, in which case it becomes
possible to bypass the sparsity assumption on $\theta^{*}$, regularization
(on the estimates of $\theta^{*}$), and other inherent challenges
in an inverse problem. As we will see in the subsequent sections,
it is quite natural that the coverage properties of our procedure
as well as the control on the Type I error do not rely on any form
of sparsity in $\theta^{*}$. In terms of the Type II error control,
we exploit that more restrictions (larger $m$) on $\theta^{*}$ in
$H_{0}$ result in fewer parameters to be determined, and if $m$
is large relative to $p-n$, it is possible for the power of our test
to not rely on any form of sparsity in $\theta^{*}$. This result
suggests that our procedure can be useful for many economic applications
where the ratio $\frac{p}{n}$ often stays constant but below $1$,
$\theta^{*}$ may not have any sparsity structure, and the null hypothesis
contains nonlinear restrictions (like the APE example given previously). 

Third, we provide a different interpretation of our approach by exploring
a pair of primal and dual optimization problems along with the so-called
Farkas's lemma (see, e.g., \cite{bt97}). In particular, the primal
(dual) problem is a feasibility (minimization) program where some
or all of entries in the target (cost) vector are ``corrupted''
with additive i.i.d. noise, e.g., Gaussian. This perspective provides
an interesting connection between our framework and Farkas' lemma
(in math programming) under uncertainty, which points to some potential
applications of our method outside traditional hypothesis testing.

\subsection{Hypothesis testing in regression models \label{subsec:Hypothesis-testing}}

We choose our test statistics in the form of
\begin{equation}
\Psi_{q}(\hat{\theta}_{\alpha}):=\left\Vert \frac{1}{n}\sum_{i=1}^{n}X_{i}\left[Y_{i}-g\left(V_{i};\hat{\theta}_{\alpha}\right)\right]\right\Vert _{q}\label{eq:stats}
\end{equation}
for some pre-specified $L-$dimensional vector of functions $f\left(V_{i}\right)=:X_{i}$.
In (\ref{eq:stats}), $\hat{\theta}_{\alpha}$ is obtained by solving
the following program\footnote{Section \ref{subsec:Alternative-formulation} provides an alternative
formulation. }: 
\begin{align}
\left(\hat{\theta}_{\alpha},\,\hat{\mu}_{\alpha}\right)\in\arg\min_{(\theta_{\alpha},\mu_{\alpha})\in\mathbb{R}^{p}\times\mathbb{R}}\mu_{\alpha}\nonumber \\
\textrm{subject to: }\left\Vert \frac{1}{n}\sum_{i=1}^{n}X_{i}\left[Y_{i}-g\left(V_{i};\theta_{\alpha}\right)\right]\right\Vert _{q}\leq r_{\alpha,q}+\mu_{\alpha},\label{eq:p2-1}\\
h(\theta_{\alpha})\in\Omega,\label{eq:5}\\
\mu_{\alpha}\geq0,\label{eq:6}
\end{align}
with $q\in\left[1,\,\infty\right]$ chosen by the users. For $1\leq q\leq\infty$,
we write $\left\Vert \Delta\right\Vert _{q}$ to mean the $l_{q}-$norm
of an $L-$dimensional vector $\Delta$, where $\left\Vert \Delta\right\Vert _{q}:=\left(\sum_{i=1}^{L}|\Delta_{i}|^{q}\right)^{1/q}$
when $1\leq q<\infty$ and $\left\Vert \Delta\right\Vert _{q}:=\max_{i=1,...,L}|\Delta_{i}|$
when $q=\infty$. Here we suppress the dependence of $\left(\hat{\theta}_{\alpha},\,\hat{\mu}_{\alpha}\right)$
in (\ref{eq:p2-1}) on $q$ for notational simplicity. The choice
for $r_{\alpha,q}$ in the first set of constraints is to be specified
in Section \ref{sec:Gaussian-Regressions}. Statistical guarantees
in this paper are stated in terms of $\left(\alpha,\,q\right)$. 

Suppose $g\left(V_{i};\theta^{*}\right)=V_{i}\theta^{*}$ in (\ref{eq:reg}).
For the special case $q=\infty$ and $X_{i}=V_{i}$, it is worth comparing
program (\ref{eq:p2-1}) with the well-known Dantzig selector:

\begin{align}
\hat{\theta}^{dan}\in\textrm{arg}\min_{\theta\in\mathbb{R}^{p}}\left\Vert \theta\right\Vert _{1}\;\;\textrm{subject to}\left\Vert \frac{1}{n}V^{T}(Y-V\theta)\right\Vert _{\infty}\leq r.\label{eq:dan}
\end{align}
Unlike (\ref{eq:dan}), (\ref{eq:p2-1}) involves a slack variable
$\mu_{\alpha}$ in the first set of constraints, (\ref{eq:5}) and
(\ref{eq:6}), as well as a different objective function (minimizing
the slack variable, instead of minimizing $\left\Vert \theta\right\Vert _{1}$).
Consequently, the resulting solution to (\ref{eq:dan}) is not constrained
to satisfy $h(\hat{\theta}^{dan})\in\Omega$, whereas $\hat{\theta}_{\alpha}$
in (\ref{eq:p2-1}) satisfies $h(\hat{\theta}_{\alpha})\in\Omega$.

A solution $\hat{\theta}_{\alpha}$ to (\ref{eq:p2-1}) may not necessarily
be unique: that is, there might be different $\hat{\theta}_{\alpha}$s
that satisfy (\ref{eq:p2-1}) while delivering the same (minimal)
objective value $\hat{\mu}_{\alpha}$. We refer to the non-negative
scalar $\mu_{\alpha}$ in (\ref{eq:p2-1}) as the ``slack'' variable
that fills the gap between $\left\Vert \frac{1}{n}\sum_{i=1}^{n}X_{i}\left[Y_{i}-g\left(V_{i};\theta^{*}\right)\right]\right\Vert _{q}$
and $\left\Vert \frac{1}{n}\sum_{i=1}^{n}X_{i}\left[Y_{i}-g\left(V_{i};\theta_{\alpha}\right)\right]\right\Vert _{q}$
where $h(\theta_{\alpha})\in\Omega$. When the null hypothesis is
true, i.e., $h(\theta^{*})\in\Omega$, the optimal value $\hat{\mu}_{\alpha}$
must be zero with probability at least $1-\alpha$. 

We derive (non-asymptotic) thresholds $r_{\alpha,q}$ such that \textit{ } 
\begin{eqnarray}
\mathbb{P}_{0}\left\{ \Psi_{q}(\hat{\theta}_{\alpha})\geq r_{\alpha,q}\right\}  & \leq & \alpha,\qquad\textrm{(Type I Error)}\label{eq:alpha}\\
\mathbb{P}_{1}\left\{ \Psi_{q}(\hat{\theta}_{\alpha})\leq r_{\alpha,q}\right\}  & \leq & \beta,\qquad\textrm{(Type II Error)}\label{eq:beta}
\end{eqnarray}
where $\mathbb{P}_{0}$ means under $H_{0}$, $\mathbb{P}_{1}$ means
under $H_{1}$ and a ``Level $\beta$ Separation Requirement'' (to
be introduced in Section \ref{subsec:Separation-requirement-and}).
Our decision rule is that if $\Psi_{q}(\hat{\theta}_{\alpha})\geq r_{\alpha,q}$,
we reject the null hypothesis $H_{0}$ in (\ref{eq:nonlinear}) at
the $1-\alpha$ confidence level. In addition to the guarantees on
the Type I and Type II errors, we also construct confidence regions
in terms of $\left(\hat{\theta}_{\alpha},\,\hat{\mu}_{\alpha}\right)$.

\subsection{Comparison with existing results}

At first glance, there seems no lack of non-asymptotic bounds on the
$l_{p}-$error (often $p\in\left[1,\,2\right]$ or $p=\infty$) of
some regularized estimator (such as the various Lasso estimators and
Dantzig selectors) concerning a sparse linear regression model (e.g.,
\cite{brt09}). However, these bounds (even in the sharpest forms)
tend to involve quite a few unknown nuisance parameters that are hard
to estimate in practice. In order to adapt the existing bounds for
the purpose of inference, knowledge on the sparsity of $\theta^{*}$
would be needed at a minimum; see the discussion in \cite{gt11}.
By contrast, the nonasymptotic thresholds in our testing procedure
do not involve any unknown parameters related to sparsity.

In particular, this paper studies nonasymptotic inference by exploiting
the concentration phenomenon, which should be distinguished from another
line of literature based on normal approximations using the Stein's
Method, for example, \cite{cck13}, \cite{h17} (also see \cite{zc17,zb17b}
whose methods are justified by the theory in \cite{cck13}). In particular,
\cite{h17} develops results for hypothesis testing in the regime
of $n\gg p$; by contrast, we also study the regime where $p$ is
comparable to or larger than $n$. In \cite{h17}, some of the results
are still only asymptotically valid and the other results (even though
nonasymptotically justified) come with probabilistic guarantees that
contain rather loose constants and dimension-dependent components. 

If $g\left(V_{i};\theta^{*}\right)=V_{i}\theta^{*}$ in (\ref{eq:reg})
and we choose $q=\infty$, $X_{i}=V_{i}$, then (\ref{eq:stats})
is reduced to 
\[
\Psi_{\infty}\left(\hat{\theta}_{\alpha}\right):=\left\Vert \frac{1}{n}\sum_{i=1}^{n}V_{i}\left(Y_{i}-V_{i}\hat{\theta}_{\alpha}\right)\right\Vert _{\infty}.
\]
This statistics shares some resemblance to the score-based correction
term in the debiased Lasso literature (see, e.g., \cite{dbz17,jm14,vbd14,zz14,zc17})
as well as the decorrelated score in \cite{nl17}. Unlike the debiased
and decorrelated procedures which require an initial (consistent)
estimator for (the sparse) $\theta^{*}$ in the correction term, our
$\hat{\theta}_{\alpha}$ here need not be consistent and is directly
used in the test statistics (requiring no further debiasing or decorrelating
step). In addition, our method is nonasymptotically valid, whereas
the aforementioned papers hinge on the asymptotic normality of the
debiased or decorrelated procedure. 

For choosing the non-asymptotic thresholds $r_{\alpha,q}$ for general
$q\in\left[1,\,\infty\right]$, our proposal exploits sharp concentration
of Lipschitz functions of Gaussian variables. The key component in
our $r_{\alpha,q}$ is data-driven and uses a Monte-Carlo approximation
to ``mimic'' the expectation that is concentrated around and automatically
captures the dependencies across coordinates. In this perspective,
our results share some similarity as those in \cite{abr10}; however,
\cite{abr10} concern inference for the mean of a random vector while
we consider inference about the coefficients ($\theta^{*}\in\mathbb{R}^{p}$)
in regression models. 

In terms of relaxing sparsity assumptions, this paper shares slight
similarity as \cite{zb17b} and \cite{zb17} although our focus and
approach are drastically different from those in \cite{zb17b} and
\cite{zb17}. First, \cite{zb17b} deal with $H_{0}:\,\theta_{j}^{*}=\theta_{j}^{0}$
for $j\in M\subseteq\left\{ 1,...,p\right\} $ and \cite{zb17} deal
with $H_{0}:\,a^{T}\theta^{*}=b^{0}$ for some prespecified $a\in\mathbb{R}^{p}$,
$b^{0}\in\mathbb{R}$. Second, parts of the analyses in \cite{zb17b}
and \cite{zb17} still rely on sparsity to some extent (albeit considerably
less compared to much of the existing literature) even in cases where
$\frac{p}{n}$ stays constant but below $1$; as a consequence, their
methods still involve regularization (on the estimates of subcomponents
of $\theta^{*}$). Finally, like much of the existing literature,
the statistical guarantees in \cite{zb17b} and \cite{zb17} are asymptotic
and require $n\rightarrow\infty$. 

\subsection{Feasibility determination under uncertainty \label{subsec:Feasibility-determination}}

The proposed framework in this paper can be motivated with a classical
notion called ``certificate of infeasibility'' in operations research.
Essentially, we can view hypothesis testing as determining a given
system of inequalities is feasible or not when some of the system
inputs are subject to uncertainty. For the exposition in this section,
let us introduce the following pair of problems:
\begin{eqnarray}
\max_{\theta\in\mathbb{R}^{p}} &  & \mathbf{0}_{p}^{T}\theta\nonumber \\
\textrm{subject to} &  & A\theta=b\label{eq:primal}\\
 &  & \theta\geq0\nonumber 
\end{eqnarray}
and 
\begin{eqnarray}
\min_{\pi\in\mathbb{R}^{d}} &  & \pi^{T}b\nonumber \\
\textrm{subject to} &  & \pi^{T}A\geq\mathbf{0}_{p}^{T}.\label{eq:dual}
\end{eqnarray}

The book \cite{bt97} refers to the $d-$dimensional column vectors
in the matrix $A\in\mathbb{R}^{d\times p}$ as the resource vectors,
$b\in\mathbb{R}^{d}$ as the target vector in (\ref{eq:primal}) and
the cost vector in (\ref{eq:dual}); $\theta\in\mathbb{R}^{p}$ and
$\pi\in\mathbb{R}^{d}$ are called the vectors of decision variables
to (\ref{eq:primal}) and (\ref{eq:dual}), respectively. Note that
program (\ref{eq:primal}) is the dual of program (\ref{eq:dual})
and vice versa. By convention in linear programming, the constraints
in (\ref{eq:primal}) are called the ``standard form'' and the $d$
rows of $A$ are assumed to be linearly independent (so we must have
$d\leq p$). Algorithms designed for linear programming can be used
to determine whether a feasible solution exists for (\ref{eq:primal}). 

The pair (\ref{eq:primal})-(\ref{eq:dual}) is associated with the
famous Farkas' lemma, which states that: (1) if (\ref{eq:primal})
is infeasible, then there exists some $\pi$ feasible for (\ref{eq:dual})
such that $\pi^{T}b<0$; (2) if (\ref{eq:primal}) is feasible, then
there cannot exist a $\pi$ feasible for (\ref{eq:dual}) such that
$\pi^{T}b<0$. This lemma has an important application in asset pricing,
in particular the so called ``absence of arbitrage'' condition that
underlies much of finance theory (see the example in \cite{bt97}).

In the classical paradigm, $b$ is assumed to be precisely known and
not subject to uncertainty. Started in the early 1970s, developing
robust approaches to solving linear programming problems where data
are subject to uncertainty has been extensively studied. Some of the
well-known papers include, for example, \cite{Ben-Tal,Bert_Sim,Soyster}.
Instead of contributing new optimization theory to this literature,
our goal here is to show an interesting link between our framework
and Farkas' lemma under uncertainty. Let us consider a special case
of uncertainty where the first $n\left(\leq d\right)$ components
in $b$ are corrupted with i.i.d. Gaussian noise, in which scenario
we only observe $Y_{i}=b_{i}+W_{i}$ with $W_{i}\sim\mathcal{N}(0,\,\sigma^{2})$
for $i=1,...,n$. We are interested in the following pairs of hypotheses:
\begin{align}
H_{0}: & \textrm{ A solution to (\ref{eq:primal}) exists}\label{eq:feasible}\\
H_{1}: & \textrm{ No solution exists for (\ref{eq:primal}),}\nonumber 
\end{align}
which is equivalent to
\begin{align}
H_{0}: & \textrm{ There cannot exist a }\pi\textrm{ feasible for (\ref{eq:dual}) such that }\pi^{T}b<0\label{eq:unbounded}\\
H_{1}: & \textrm{ There exists a }\pi\textrm{ feasible for (\ref{eq:dual}) such that }\pi^{T}b<0.\nonumber 
\end{align}

To test the hypotheses above, we can adapt (\ref{eq:p2-1}) and solve
the following program: 
\begin{align}
\left(\hat{\theta}_{\alpha},\,\hat{\mu}_{\alpha}\right)\in\arg\min_{(\theta_{\alpha},\mu_{\alpha})\in\mathbb{R}^{p}\times\mathbb{R}}\mu_{\alpha}\nonumber \\
\textrm{subject to: }\left\Vert \frac{1}{n}\sum_{i=1}^{n}X_{i}\left(Y_{i}-V_{i}\theta_{\alpha}\right)\right\Vert _{q}\leq & r_{\alpha,q}+\mu_{\alpha},\label{eq:lp}\\
A_{i}\theta_{\alpha}=b_{i}, & i=n+1,...,d,\nonumber \\
\theta_{\alpha}\geq0,\nonumber \\
\mu_{\alpha}\geq0,\nonumber 
\end{align}
for some pre-specified $L-$dimensional vector of functions $f\left(V_{i}\right)=:X_{i}$,
where $V_{i}=A_{i}$ for $i=1,...,n$ and $A_{i}$ is the $i$th row
of $A$. Note that program (\ref{eq:lp}) is simply a special instance
of (\ref{eq:p2-1}), where the constraints $h(\theta_{\alpha})\in\Omega$
correspond to 
\[
\left\{ \theta_{\alpha}\geq0,\,A_{i}\theta_{\alpha}=b_{i},\,i=n+1,...,d\right\} .
\]
\footnote{If $\left\{ \theta_{\alpha}\geq0,\,A_{i}\theta_{\alpha}=b_{i},\,i=n+1,...,d\right\} =\emptyset$,
then we can immediately conclude that no solution exists for program
(\ref{eq:primal}). As a result, like in (\ref{eq:nonlinear}), we
may assume $\left\{ \theta_{\alpha}\geq0,\,A_{i}\theta_{\alpha}=b_{i},\,i=n+1,...,d\right\} \neq\emptyset$.}Letting 
\begin{equation}
\Psi_{q}(\hat{\theta}_{\alpha}):=\left\Vert \frac{1}{n}\sum_{i=1}^{n}X_{i}\left(Y_{i}-V_{i}\hat{\theta}_{\alpha}\right)\right\Vert _{q},\label{eq:stats-1}
\end{equation}
if $\Psi_{q}(\hat{\theta}_{\alpha})\geq r_{\alpha,q}$, we then reject
(\ref{eq:feasible}) and (\ref{eq:unbounded}) at the $1-\alpha$
confidence level. As in Section \ref{subsec:Hypothesis-testing},
(\ref{eq:alpha}) and (\ref{eq:beta}) also hold for (\ref{eq:stats-1}).
This result can be viewed as the Farkas' lemma extended to cases where
the target vector $b$ in (\ref{eq:primal}) (respectively, the cost
vector in (\ref{eq:dual})) is corrupted with i.i.d. noise.

\subsection{Organization of this paper}

Section 2 establishes nonasymptotic control on the Type I and Type
II errors for our method. We provide numerical evidence through a
simulation study in Section 3. Section 4 provides additional extensions,
which include some nonasymptotic justifications for inference in statistical
models that involve non-Gaussian responses. All technical details
are collected in the supplementary materials.

\section{Main results \label{sec:Gaussian-Regressions}}

For the regression model (\ref{eq:reg}), we begin with the scenario
where $\sigma^{2}$ is known, and then consider the scenario where
$\sigma^{2}$ is unknown. Throughout this section, we use $\mathbb{E}_{W}\left[\cdot\right]$
to denote the expectation over $W$ only, conditioning on $V$.

Recalling the pre-specified $L-$dimensional vector of functions $f\left(V_{i}\right)=:X_{i}$
in Section \ref{subsec:Hypothesis-testing}, our first result establishes
an ``ideal'' confidence region for 
\[
\left\Vert \frac{1}{n}\sum_{i=1}^{n}X_{i}\left[g\left(V_{i};\theta^{*}\right)-g\left(V_{i};\hat{\theta}_{\alpha}^{*}\right)\right]\right\Vert _{q}.
\]
The ``theoretical'' optimal solution above, $\hat{\theta}_{\alpha}^{*}$,
is obtained by setting $r_{\alpha,q}$ in (\ref{eq:p2-1}) to $\mathbb{E}_{W}\left[\left\Vert \frac{1}{n}X^{T}W\right\Vert _{q}\right]$
plus a deviation. In practice, $\mathbb{E}_{W}\left[\left\Vert \frac{1}{n}X^{T}W\right\Vert _{q}\right]$
may be bounded with its Monte Carlo approximation and a ``small''
deviation term. This approach results in a ``practical'' optimal
solution, $\hat{\theta}_{\alpha}$, which can then be used to construct
test statistics and a ``practical'' confidence region.

To state the first result, we introduce the following notation (which
will appear in many places throughout this paper): 
\begin{eqnarray*}
\left\Vert \sqrt{\frac{1}{n}\sum_{i=1}^{n}X_{i}^{2}}\right\Vert _{q} & = & \sqrt[q]{\sum_{j=1}^{L}\left(\sqrt{\frac{1}{n}\sum_{i=1}^{n}X_{ij}^{2}}\right)^{q}},\qquad q\in[1,\,\infty)\\
\left\Vert \sqrt{\frac{1}{n}\sum_{i=1}^{n}X_{i}^{2}}\right\Vert _{q} & = & \max_{j\in\left\{ 1,...,L\right\} }\sqrt{\frac{1}{n}\sum_{i=1}^{n}X_{ij}^{2}},\qquad q=\infty.
\end{eqnarray*}
We can assume $X$ is normalized such that $\frac{1}{n}\sum_{i=1}^{n}X_{ij}^{2}=1$
for all $j=1,...,L$, in which case $\left\Vert \sqrt{\frac{1}{n}\sum_{i=1}^{n}X_{i}^{2}}\right\Vert _{q}=L^{\frac{1}{q}}$
for $q\in[1,\,\infty)$ and $\left\Vert \sqrt{\frac{1}{n}\sum_{i=1}^{n}X_{i}^{2}}\right\Vert _{q}=1$
for $q=\infty$.\\
\textbf{\textit{}}\\
\textbf{\textit{Proposition 2.1}}\textit{. Assume (\ref{eq:reg})
where $W\sim\mathcal{N}(\mathbf{0}_{n},\,\sigma^{2}\boldsymbol{I}_{n})$
and is independent of $V$. Then for any $q\in\left[1,\,\infty\right]$,
we have 
\begin{equation}
\mathbb{P}\left\{ \left\Vert \frac{1}{n}X^{T}W\right\Vert _{q}\geq\mathbb{E}_{W}\left[\left\Vert \frac{1}{n}X^{T}W\right\Vert _{q}\right]+t\right\} \leq\exp\left(\frac{-nt^{2}}{2\sigma^{2}\left\Vert \sqrt{\frac{1}{n}\sum_{i=1}^{n}X_{i}^{2}}\right\Vert _{q}^{2}}\right).\label{eq:prop2-1}
\end{equation}
Moreover, for $\alpha\in(0,\,1)$, let }

\textit{
\begin{equation}
r_{\alpha,q}=r_{\alpha,q}^{*}:=\mathbb{E}_{W}\left[\left\Vert \frac{1}{n}X^{T}W\right\Vert _{q}\right]+\sigma\left\Vert \sqrt{\frac{1}{n}\sum_{i=1}^{n}X_{i}^{2}}\right\Vert _{q}\sqrt{\frac{2}{n}\log\frac{1}{\alpha}}\label{eq:8}
\end{equation}
in (\ref{eq:p2-1}). Then, an optimal solution $\left(\hat{\theta}_{\alpha}^{*},\,\hat{\mu}_{\alpha}^{*}\right)$
to (\ref{eq:p2-1}) must satisfy }

\textit{
\begin{eqnarray}
\left\Vert \frac{1}{n}\sum_{i=1}^{n}X_{i}\left[g\left(V_{i};\theta^{*}\right)-g\left(V_{i};\hat{\theta}_{\alpha}^{*}\right)\right]\right\Vert _{q} & \geq & \hat{\mu}_{\alpha}^{*},\label{eq:low2}\\
\left\Vert \frac{1}{n}\sum_{i=1}^{n}X_{i}\left[g\left(V_{i};\theta^{*}\right)-g\left(V_{i};\hat{\theta}_{\alpha}^{*}\right)\right]\right\Vert _{q} & \leq & 2r_{\alpha,q}^{*}+\hat{\mu}_{\alpha}^{*},\label{eq:6-1}
\end{eqnarray}
with probability at least $1-\alpha$. }

\subsection{Hypothesis testing }

For the moment, suppose we set $r_{\alpha,q}=r_{\alpha,q}^{*}$ in
(\ref{eq:p2-1})\textcolor{red}{{} }according to (\ref{eq:8}). Under
$H_{0}$, $\left(\theta^{*},\,0\right)$ is an optimal solution to
(\ref{eq:p2-1}). Consequently, given the test statistics (\ref{eq:stats})
and a chosen $\alpha\in\left(0,\,1\right)$, any optimal solution
to (\ref{eq:p2-1}) must satisfy 
\begin{equation}
\mathbb{P}_{0}\left\{ \Psi_{q}(\hat{\theta}_{\alpha}^{*})\geq r_{\alpha,q}^{*}\right\} \leq\alpha\label{eq:7-1}
\end{equation}
where $\mathbb{P}_{0}$ means under $H_{0}$.

The claim in (\ref{eq:7-1}) suggests a test (with level $\alpha$)
based on the statistics $\Psi_{q}(\hat{\theta}_{\alpha}^{*})$ and
an ``ideal'' critical value, $r_{\alpha,q}^{*}$, given in (\ref{eq:8}).
When $W\sim\mathcal{N}(\mathbf{0}_{n},\,\sigma^{2}\boldsymbol{I}_{n})$
and $\sigma^{2}$ is known, the first term $\mathbb{E}_{W}\left[\left\Vert \frac{1}{n}X^{T}W\right\Vert _{q}\right]$
in $r_{\alpha,q}^{*}$ can be approximated by Monte-Carlo as follows.
Let $Z\in\mathbb{R}^{n\times R}$ (independent of $W$ and $V$) be
a matrix consisting of independent entries randomly drawn from $\mathcal{N}(0,\,1)$
and the $r$th column of $Z$ is denoted by $Z_{r}$. By (\ref{eq:subn})
and (\ref{eq:mgfg}) in Section \ref{subsec:Preliminary}, note that
$\sigma R^{-1}\sum_{r=1}^{R}\left\Vert \frac{1}{n}X^{T}Z_{r}\right\Vert _{q}-\mathbb{E}_{W}\left[\left\Vert \frac{1}{n}X^{T}W\right\Vert _{q}\right]$
is sub-Gaussian with parameter at most $(nR)^{-1/2}\sigma\left\Vert \sqrt{\frac{1}{n}\sum_{i=1}^{n}X_{i}^{2}}\right\Vert _{q}$.
Consequently, (\ref{eq:lower}) yields the following concentration
\begin{equation}
\mathbb{P}\left\{ \mathbb{E}_{W}\left[\left\Vert \frac{1}{n}X^{T}W\right\Vert _{q}\right]\geq\frac{\sigma}{R}\sum_{r=1}^{R}\left\Vert \frac{1}{n}X^{T}Z_{r}\right\Vert _{q}+t\right\} \leq\exp\left(\frac{-nRt^{2}}{2\sigma^{2}\left\Vert \sqrt{\frac{1}{n}\sum_{i=1}^{n}X_{i}^{2}}\right\Vert _{q}^{2}}\right).\label{eq:c2}
\end{equation}
Combining (\ref{eq:prop2-1}) and (\ref{eq:c2}) yields 
\begin{align}
 & \mathbb{P}\left\{ \left\Vert \frac{1}{n}X^{T}W\right\Vert _{q}\geq\frac{\sigma}{R}\sum_{r=1}^{R}\left\Vert \frac{1}{n}X^{T}Z_{r}\right\Vert _{q}+t_{1}+t_{2}\right\} \nonumber \\
\leq & \exp\left(\frac{-nt_{1}^{2}}{2\sigma^{2}\left\Vert \sqrt{\frac{1}{n}\sum_{i=1}^{n}X_{i}^{2}}\right\Vert _{q}^{2}}\right)+\exp\left(\frac{-nRt_{2}^{2}}{2\sigma^{2}\left\Vert \sqrt{\frac{1}{n}\sum_{i=1}^{n}X_{i}^{2}}\right\Vert _{q}^{2}}\right).\label{eq:9}
\end{align}

\subsubsection{Construction of critical values ($r_{\alpha,q}$) and Type I error}

For some chosen $\alpha_{1},\,\alpha_{2}>0$ such that $\alpha_{1}+\alpha_{2}=\alpha\in(0,\,1)$,
we let in (\ref{eq:9}), 
\begin{eqnarray}
t_{1} & = & \sigma\left\Vert \sqrt{\frac{1}{n}\sum_{i=1}^{n}X_{i}^{2}}\right\Vert _{q}\sqrt{\frac{2}{n}\log\frac{1}{\alpha_{1}}}:=\tau_{\alpha_{1},q},\label{eq:t1}\\
t_{2} & = & \sigma\left\Vert \sqrt{\frac{1}{n}\sum_{i=1}^{n}X_{i}^{2}}\right\Vert _{q}\sqrt{\frac{2}{nR}\log\frac{1}{\alpha_{2}}}:=\sqrt{\frac{1}{R}}\tau_{\alpha_{2},q}.\nonumber 
\end{eqnarray}
Based on (\ref{eq:9}) along with the choices of $t_{1}$ and $t_{2}$
above, we set in (\ref{eq:p2-1}),
\begin{equation}
r_{\alpha,q}=\frac{\sigma}{R}\sum_{r=1}^{R}\left\Vert \frac{1}{n}X^{T}Z_{r}\right\Vert _{q}+\tau_{\alpha_{1},q}+\sqrt{\frac{1}{R}}\tau_{\alpha_{2},q}.\label{eq:3}
\end{equation}
Note that we can draw as many columns in $Z$ as we want, to make
$\sqrt{\frac{1}{R}}\tau_{\alpha_{2},q}$ in (\ref{eq:3}) small; for
a given $\alpha$, we can let $\alpha_{2}$ be smaller than $\alpha_{1}$
because of the additional ``$\sqrt{\frac{1}{R}}$''. 

Under $H_{0}$, $\left(\theta^{*},\,0\right)$ is an optimal solution
to (\ref{eq:p2-1}) with $r_{\alpha,q}$ specified in (\ref{eq:3}).
Consequently, a (practical) optimal solution to (\ref{eq:p2-1}) must
satisfy \textit{ 
\begin{equation}
\mathbb{P}_{0}\left\{ \Psi_{q}(\hat{\theta}_{\alpha})\geq r_{\alpha,q}\right\} \leq\alpha\qquad\textrm{(Type I Error)}.\label{eq:13}
\end{equation}
}\textbf{Remarks}. In terms of control on the Type I error, the $l_{q}-$norm
in (\ref{eq:stats}) and (\ref{eq:p2-1}) can be generalized to the
function $\zeta_{q}:\,\mathbb{R}^{L}\mapsto\mathbb{R}$ that satisfies: 
\begin{itemize}
\item for all $z\in\mathbb{R}^{L}$ and $a\in\mathbb{R}^{+}$, $\zeta_{q}\left(az\right)=a\zeta_{q}\left(z\right)$,
\item for all $z,\,z^{'}\in\mathbb{R}^{L}$, $\zeta_{q}\left(z+z^{'}\right)\leq\zeta_{q}\left(z\right)+\zeta_{q}(z^{'})$,
\item for all $z\in\mathbb{R}^{L}$, $\left|\zeta_{q}\left(z\right)\right|\leq\left\Vert z\right\Vert _{q}$
for $q\in\left[1,\,\infty\right]$. 
\end{itemize}
This generalization is inspired by the construction in \cite{abr10}.
We simply let 
\[
r_{\alpha,q}=\frac{\sigma}{R}\sum_{r=1}^{R}\zeta_{q}\left(\frac{1}{n}X^{T}Z_{r}\right)+\tau_{\alpha_{1},q}+\sqrt{\frac{1}{R}}\tau_{\alpha_{2},q},
\]
and obtain 
\[
\mathbb{P}_{0}\left\{ \zeta_{q}\left(\frac{1}{n}\sum_{i=1}^{n}X_{i}\left[Y_{i}-g\left(V_{i};\hat{\theta}_{\alpha}\right)\right]\right)\geq r_{\alpha,q}\right\} \leq\alpha\qquad\textrm{(Type I Error)},
\]
where $\hat{\theta}_{\alpha}$ is a solution to (\ref{eq:p2-1}) with
the $l_{q}-$norm in the first set of constraints replaced by $\zeta_{q}$.
Given $\zeta_{q}$ is subadditive and bounded by the $l_{q}-$norm,
the result above follows from the simple fact that
\begin{eqnarray*}
\left|\zeta_{q}\left(\frac{1}{n}X^{T}W\right)-\zeta_{q}\left(\frac{1}{n}X^{T}W^{'}\right)\right| & \leq & \left\Vert \frac{1}{n}X^{T}\left(W-W^{'}\right)\right\Vert _{q}\\
 & \leq & \frac{1}{\sqrt{n}}\left\Vert \sqrt{\frac{1}{n}\sum_{i=1}^{n}X_{i}^{2}}\right\Vert _{q}\left\Vert W-W^{'}\right\Vert _{2}.
\end{eqnarray*}
Consequently, we can establish bounds that are identical to (\ref{eq:prop2-1}),
(\ref{eq:c2}), (\ref{eq:9}) in terms of $\zeta_{q}\left(\frac{1}{n}X^{T}W\right)$,
$\mathbb{E}_{W}\left[\zeta_{q}\left(\frac{1}{n}X^{T}W\right)\right]$,
$\frac{\sigma}{R}\sum_{r=1}^{R}\zeta_{q}\left(\frac{1}{n}X^{T}Z_{r}\right)$,
and then follow the same argument as what is used to show (\ref{eq:13}).

\subsubsection{Practical confidence regions }

Let $\left(\hat{\theta}_{\alpha},\,\hat{\mu}_{\alpha}\right)$ be
an optimal solution to (\ref{eq:p2-1}) with $r_{\alpha,q}$ specified
in (\ref{eq:3}). We have 
\begin{eqnarray}
 &  & \left\Vert \frac{1}{n}\sum_{i=1}^{n}X_{i}\left[g\left(V_{i};\theta^{*}\right)-g\left(V_{i};\hat{\theta}_{\alpha}\right)\right]\right\Vert _{q}\nonumber \\
 & \leq & \left\Vert \frac{1}{n}\sum_{i=1}^{n}X_{i}\left[Y_{i}-g\left(V_{i};\hat{\theta}_{\alpha}\right)\right]\right\Vert _{q}+\left\Vert \frac{1}{n}X^{T}W\right\Vert _{q}\nonumber \\
 & \leq & \frac{2\sigma}{R}\sum_{r=1}^{R}\left\Vert \frac{1}{n}X^{T}Z_{r}\right\Vert _{q}+2\tau_{\alpha_{1},q}+2\sqrt{\frac{1}{R}}\tau_{\alpha_{2},q}+\hat{\mu}_{\alpha}\label{eq:7b}
\end{eqnarray}
and 
\begin{equation}
\left\Vert \frac{1}{n}\sum_{i=1}^{n}X_{i}\left[g\left(V_{i};\theta^{*}\right)-g\left(V_{i};\hat{\theta}_{\alpha}\right)\right]\right\Vert _{q}\geq\hat{\mu}_{\alpha}\label{eq:7b-low}
\end{equation}
with probability at least $1-\alpha$. The argument for (\ref{eq:7b-low})
is identical to what is used to show (\ref{eq:low2}) (see more details
in Section \ref{subsec:Proof-of-Propositions 2.1 and A2.1}). As we
have pointed out in the introduction, there might be different $\hat{\theta}_{\alpha}$s
that solve (\ref{eq:p2-1}) while producing the same (minimal) objective
value $\hat{\mu}_{\alpha}$. Consequently, there may be more than
one confidence interval in the form of (\ref{eq:7b}) and (\ref{eq:7b-low}).
In view of these intervals, the length of them is naturally
\begin{equation}
CI-Length=\frac{2\sigma}{R}\sum_{r=1}^{R}\left\Vert \frac{1}{n}X^{T}Z_{r}\right\Vert _{q}+2\tau_{\alpha_{1},q}+2\sqrt{\frac{1}{R}}\tau_{\alpha_{2},q}.\label{eq:length}
\end{equation}

If $\mathbb{E}_{W}\left[\left\Vert \frac{1}{n}X^{T}W\right\Vert _{q}\right]$
can be known exactly and we were able to set $r_{\alpha_{1},q}=r_{\alpha_{1},q}^{*}$
in (\ref{eq:p2-1}) as in Proposition 2.1, then any resulting optimal
solution $\left(\hat{\theta}_{\alpha}^{*},\,\hat{\mu}_{\alpha}^{*}\right)$
to (\ref{eq:p2-1}) should satisfy 
\begin{equation}
\left\Vert \frac{1}{n}\sum_{i=1}^{n}X_{i}\left[g\left(V_{i};\theta^{*}\right)-g\left(V_{i};\hat{\theta}_{\alpha}^{*}\right)\right]\right\Vert _{q}\leq2\mathbb{E}_{W}\left[\left\Vert \frac{1}{n}X^{T}W\right\Vert _{q}\right]+2\tau_{\alpha_{1},q}+\hat{\mu}_{\alpha}^{*}\label{eq:14b}
\end{equation}
with probability at least $1-\alpha_{1}$. Comparing (\ref{eq:7b})
with (\ref{eq:14b}), note that the difference in the right hand sides
is 
\[
2\left(\frac{\sigma}{R}\sum_{r=1}^{R}\left\Vert \frac{1}{n}X^{T}Z_{r}\right\Vert _{q}-\mathbb{E}_{W}\left[\left\Vert \frac{1}{n}X^{T}W\right\Vert _{q}\right]\right)+2\sqrt{\frac{1}{R}}\tau_{\alpha_{2},q},
\]
which can be made arbitrarily small with a large number of random
draws in the Monte-Carlo approximation. Because of such an approximation,
the probabilistic guarantee for (\ref{eq:7b}) is bounded from below
by $1-\alpha$ instead of $1-\alpha_{1}$.

Given the statistics $\Psi_{q}(\hat{\theta}_{\alpha})$ in (\ref{eq:stats})
based on (a practical) $\hat{\theta}_{\alpha}$ and the critical value
$r_{\alpha,q}$ defined in (\ref{eq:3}), we have constructed a test
with level $\alpha$ as shown in (\ref{eq:13}). For a chosen $\beta\in(0,\,1)$,
when can this test correctly detect an alternative with probability
at least $1-\beta$? To answer this question, we introduce the ``Separation
Requirement'' in the following section.

\subsubsection{Separation requirement and Type II error\label{subsec:Separation-requirement-and}}

Letting $\Theta_{0}:=\left\{ \theta\in\mathbb{R}^{p}\,:\,h(\theta)\in\Omega\right\} $,
we choose $\beta_{1},\,\beta_{2}>0$ such that $\beta_{1}+\beta_{2}=\beta\in(0,\,1)$,
and assume 
\begin{equation}
\inf_{\theta\in\Theta_{0}}\left\Vert \frac{1}{n}\sum_{i=1}^{n}X_{i}\left[g\left(V_{i};\theta^{*}\right)-g\left(V_{i};\theta\right)\right]\right\Vert _{q}\geq\delta_{\alpha,\beta,q}\label{eq:ms}
\end{equation}
with 
\begin{equation}
\delta_{\alpha,\beta,q}=2\mathbb{E}_{W}\left[\left\Vert \frac{1}{n}X^{T}W\right\Vert _{q}\right]+\tau_{\alpha_{1},q}+\sqrt{\frac{1}{R}}\tau_{\alpha_{2},q}+\sqrt{\frac{1}{R}}\tau_{\beta_{1},q}+\tau_{\beta_{2},q}\label{eq:ms1}
\end{equation}
for the prespecified $\alpha_{1},\,\alpha_{2}>0$ (used in (\ref{eq:3}))
such that $\alpha_{1}+\alpha_{2}=\alpha\in(0,\,1)$. We will refer
to (\ref{eq:ms}) as the ``Separation Requirement'' at the level
$\beta$. 

Our next result concerns the Type II error of the test based on $\Psi_{q}(\hat{\theta}_{\alpha})$
in (\ref{eq:stats}) and $r_{\alpha,q}$ defined in (\ref{eq:3}).
For completeness, we also include the result for the Type I error.\\
 \\
\textbf{\textit{Theorem 2.1}}\textit{. Assume (\ref{eq:reg}) where
$W\sim\mathcal{N}(\mathbf{0}_{n},\,\sigma^{2}\boldsymbol{I}_{n})$
and is independent of $V$. For chosen $\alpha_{1},\,\alpha_{2}>0$
such that $\alpha_{1}+\alpha_{2}=\alpha\in(0,\,1)$, consider the
statistics $\Psi_{q}(\hat{\theta}_{\alpha})$ based on $\hat{\theta}_{\alpha}$
and the critical value $r_{\alpha,q}$ defined in (\ref{eq:3}). For
any $q\in\left[1,\,\infty\right]$, we have 
\begin{equation}
\mathbb{P}_{0}\left\{ \Psi_{q}(\hat{\theta}_{\alpha})\geq r_{\alpha,q}\right\} \leq\alpha,\qquad\textrm{(Type I Error)}\label{eq:13-2}
\end{equation}
where $\mathbb{P}_{0}$ means under $H_{0}$. For the same $r_{\alpha,q}$
used in (\ref{eq:13-2}) and chosen $\beta_{1},\,\beta_{2}>0$ such
that $\beta_{1}+\beta_{2}=\beta\in(0,\,1)$, if $h(\theta^{*})\notin\Omega$
and (\ref{eq:ms}) is satisfied, we have
\begin{equation}
\mathbb{P}_{1}\left\{ \Psi_{q}(\hat{\theta}_{\alpha})\leq r_{\alpha,q}\right\} \leq\beta,\qquad\textrm{(Type II Error)}\label{eq:13-1}
\end{equation}
where $\mathbb{P}_{1}$ means under $H_{1}$ and (\ref{eq:ms}). }

\subsubsection{Discussions of the results\label{subsec:Discussions}}

Some observations can be made from the results we have established
so far. First, our guarantees do not rely on any form of sparsity
in $\theta^{*}$, well behaved $\left\Vert \hat{\theta}_{\alpha}-\theta^{*}\right\Vert _{2}$
or $\sqrt{\frac{\sum_{i=1}^{n}\left[g\left(V_{i};\theta^{*}\right)-g\left(V_{i};\hat{\theta}\right)\right]^{2}}{n}}$. 

Second, the number of restrictions (i.e., $m$) in $H_{0}$ plays
a significant role in our separation requirement. Suppose $\Theta_{0}\neq\emptyset$.
For (\ref{eq:ms}) to hold, it is almost necessary that $m>p-n$.
If $p\leq n$, this ``necessary'' condition is satisfied for any
$m>0$. If $p>n$ but $m\leq p-n$, we can always find a solution
$\hat{\theta}$ such that $h(\hat{\theta})\in\Omega$ and $g\left(V_{i};\hat{\theta}\right)=Y_{i}$
for all $i$. Consequently, we have 
\begin{align}
 & \left\Vert \frac{1}{n}\sum_{i=1}^{n}X_{i}\left[g\left(V_{i};\theta^{*}\right)-g\left(V_{i};\hat{\theta}\right)\right]\right\Vert _{q}\nonumber \\
= & \left\Vert \frac{1}{n}\sum_{i=1}^{n}X_{i}\left[Y_{i}-g\left(V_{i};\hat{\theta}\right)\right]-\frac{1}{n}\sum_{i=1}^{n}X_{i}\left[Y_{i}-g\left(V_{i};\theta^{*}\right)\right]\right\Vert _{q}=\left\Vert \frac{1}{n}X^{T}W\right\Vert _{q}.\label{eq:small}
\end{align}
By (\ref{eq:prop2-1}), 
\[
\mathbb{P}\left(\left\Vert \frac{1}{n}X^{T}W\right\Vert _{q}\leq\mathbb{E}_{W}\left[\left\Vert \frac{1}{n}X^{T}W\right\Vert _{q}\right]+\tau_{\alpha,q}\right)\geq1-\alpha.
\]
Since $\delta_{\alpha,\beta,q}>\mathbb{E}_{W}\left[\left\Vert \frac{1}{n}X^{T}W\right\Vert _{q}\right]+\tau_{\alpha,q}$,
for a reasonably small Type I error $\alpha>0$, $\mathbb{P}\left(\left\Vert \frac{1}{n}\sum_{i=1}^{n}X_{i}\left[g\left(V_{i};\theta^{*}\right)-g\left(V_{i};\hat{\theta}\right)\right]\right\Vert _{q}\geq\delta_{\alpha,\beta,q}\right)$
is small and our procedure is unlikely to detect the alternatives.
As $m$ gets larger relative to $p-n$, it becomes easier for (\ref{eq:ms})
to be satisfied and our procedure to detect the alternative. Note
that $\left(\hat{\theta},\,0\right)$ also solves (\ref{eq:p2-1})
with probability $1$ for any $r_{\alpha,q}\geq0$, and clearly the
control on the Type I error remains valid. 

It is instructive to consider the case $q=\infty$. Whenever $p-n<m<p$,
for (\ref{eq:ms}) to hold, it is almost necessary that $L+m>p$ (recalling
$L$ is the dimension of $X_{i}$). Suppose $L+m\leq p$. We can always
find a solution $\hat{\theta}$ such that $\frac{1}{n}\sum_{i=1}^{n}X_{i}\left[Y_{i}-g\left(V_{i};\hat{\theta}\right)\right]=\mathbf{0}_{L}$
and $h(\hat{\theta})\in\Omega$. Once again, this fact leads to (\ref{eq:small})
and a similar consequence. 

Letting $\hat{\theta}_{\alpha}$ be a solution that gives the optimal
value $\hat{\mu}_{\alpha}$ and the $L-$dimensional vector
\begin{eqnarray*}
\tilde{\mu}_{\alpha} & = & \frac{1}{n}\sum_{i=1}^{n}X_{i}\left[Y_{i}-g\left(V_{i};\hat{\theta}_{\alpha}\right)\right]-\frac{1}{n}\sum_{i=1}^{n}X_{i}\left[Y_{i}-g\left(V_{i};\theta^{*}\right)\right]\\
 & = & \frac{1}{n}\sum_{i=1}^{n}X_{i}\left[g\left(V_{i};\theta^{*}\right)-g\left(V_{i};\hat{\theta}_{\alpha}\right)\right],
\end{eqnarray*}
Section \ref{subsec:Proof-of-Propositions 2.1 and A2.1} then implies
that $\left(\hat{\theta}_{\alpha},\,\left\Vert \tilde{\mu}_{\alpha}\right\Vert _{\infty}\right)$
is a feasible solution to (\ref{eq:p2-1}) with probability at least
$1-\alpha$. Consequently, the optimality of $\left(\hat{\theta}_{\alpha},\,\hat{\mu}_{\alpha}\right)$
implies
\[
\left\Vert \tilde{\mu}_{\alpha}\right\Vert _{\infty}=\left\Vert \frac{1}{n}\sum_{i=1}^{n}X_{i}\left[g\left(V_{i};\theta^{*}\right)-g\left(V_{i};\hat{\theta}_{\alpha}\right)\right]\right\Vert _{\infty}\geq\hat{\mu}_{\alpha}.
\]
In summary, as $\hat{\mu}_{\alpha}$ increases, the actual separation
$\left\Vert \tilde{\mu}_{\alpha}\right\Vert _{\infty}$ will never
decrease and it becomes easier for our procedure to detect the alternatives.

This fact has an implication on the choices of $X_{i}=f\left(V_{i}\right)$,
which can be quite flexible in (\ref{eq:p2-1}) given $W_{i}$ is
independent of $V_{i}$. Note that including more components in $X_{i}$
may or may not increase the optimal value $\hat{\mu}_{\alpha}$ as
larger $L$ introduces more constraints in (\ref{eq:p2-1}) but at
the same time increases $r_{\alpha,\infty}$. However, the latter
effect might be dominated by the former when the increase in $L$
is moderate: as we will see in (\ref{eq:lead1}), one can show that
$\mathbb{E}_{W}\left[\left\Vert \frac{1}{n}X^{T}W\right\Vert _{\infty}\right]\precsim\sqrt{\frac{\log L}{n}}$,
and by (\ref{eq:c2_u}) of Section \ref{subsec:Proof-of-Theorem 2.1},
\begin{equation}
r_{\alpha,\infty}\precsim\sqrt{\frac{\log L}{n}}+\sigma\sqrt{\frac{2}{n}\log\frac{1}{\alpha_{1}}}+\sigma\sqrt{\frac{2}{nR}\log\frac{1}{\alpha_{2}}}\label{eq:33}
\end{equation}
with high probability, where we have imposed the normalization on
$X$ such that $\frac{1}{n}\sum_{i=1}^{n}X_{ij}^{2}=1$ for all $j=1,...,L$.
Because of the ``$\log L$'' factor in (\ref{eq:33}), under the
alternative hypothesis, a somewhat larger $L$ may shrink the feasible
region for (\ref{eq:p2-1}) and increase the optimal value $\hat{\mu}_{\alpha}$
as well as the actual separation $\left\Vert \frac{1}{n}\sum_{i=1}^{n}X_{i}\left[g\left(V_{i};\theta^{*}\right)-g\left(V_{i};\hat{\theta}_{\alpha}\right)\right]\right\Vert _{\infty}$. 

\subsection{The union bound alternative\label{subsec:The-Bonferroni-alternative}}

As an alternative, the union bound approach can also be used to construct
a testing procedure. In particular, we can solve (\ref{eq:p2-1})
with $q=\infty$ and 
\begin{equation}
r_{\alpha,\infty}=\sqrt{\max_{j\in\left\{ 1,...,L\right\} }\frac{2\sigma^{2}}{n}\sum_{i=1}^{n}X_{ij}^{2}}\sqrt{\frac{1}{n}\log\frac{2L}{\alpha}}.\label{eq:bon}
\end{equation}
Consequently, the separation distance in (\ref{eq:ms}) that allows
us to correctly detect an alternative with probability at least $1-\beta$
takes the form 
\begin{eqnarray}
\delta_{\alpha,\beta,\infty} & = & r_{\alpha,\infty}+r_{\beta,\infty}\nonumber \\
 & = & \sqrt{\max_{j\in\left\{ 1,...,L\right\} }\frac{2\sigma^{2}}{n}\sum_{i=1}^{n}X_{ij}^{2}}\left(\sqrt{\frac{1}{n}\log\frac{2L}{\alpha}}+\sqrt{\frac{1}{n}\log\frac{2L}{\beta}}\right).\label{eq:ms-bon}
\end{eqnarray}
Given (\ref{eq:bon}) and (\ref{eq:ms-bon}), we can apply the same
argument for Theorem 2.1 and arrive at (\ref{eq:13-2}) and (\ref{eq:13-1})
(with $q=\infty$ in both). In contrast to our previous concentration
approach, (\ref{eq:bon}) is derived from a simple union bound on
$\left\Vert \frac{1}{n}X^{T}W\right\Vert _{\infty}$; as a consequence,
the resulting threshold $r_{\alpha,\infty}$ does not capture the
dependencies between the coordinates. 

We observe from (\ref{eq:t1}), (\ref{eq:3}) and (\ref{eq:ms1})
that the quantities taking the form of $\sqrt{\log\frac{1}{\varsigma}}$
are dimension free in the concentration approach. Instead, the leading
term $\mathbb{E}_{W}\left[\left\Vert \frac{1}{n}X^{T}W\right\Vert _{q}\right]$
in (\ref{eq:prop2-1}) and (\ref{eq:ms1}) reflects the ``dimension
complexity'' and automatically takes into consideration the dependencies
between the coordinates. This result is a direct consequence of the
concentration phenomenon in Lipschitz functions of Gaussians. Suppose
$X$ is normalized such that $\frac{1}{n}\sum_{i=1}^{n}X_{ij}^{2}=1$
for all $j=1,...,L$, and $\tau_{\varsigma,\infty}=\sigma\sqrt{\frac{2}{n}\log\frac{1}{\varsigma}}$.
Take $q=\infty$, $W\sim\mathcal{N}\left(\mathbf{0}_{n},\,\boldsymbol{I}_{n}\right)$
and consider the extreme example where $X$ consists of $L$ copies
of the same column $X_{0}$. Then, we have 
\[
\mathbb{E}_{W}\left[\left\Vert \frac{1}{n}X^{T}W\right\Vert _{\infty}\right]=\sqrt{\frac{2}{\pi}}\frac{1}{\sqrt{n}}
\]
and (\ref{eq:ms1}) becomes 
\begin{equation}
\delta_{\alpha,\beta,\infty}=2\sqrt{\frac{2}{\pi}}\frac{1}{\sqrt{n}}+\tau_{\alpha_{1},\infty}+\sqrt{\frac{1}{R}}\tau_{\alpha_{2},\infty}+\sqrt{\frac{1}{R}}\tau_{\beta_{1},\infty}+\tau_{\beta_{2},\infty},\label{eq:free}
\end{equation}
which involves no dimension complexity (as desired). 

Beyond the extreme example, more generally for $q=\infty$ and $W\sim\mathcal{N}\left(\mathbf{0}_{n},\,\boldsymbol{I}_{n}\right)$
(without much loss of generality by assuming $\sigma=1$), Section
\ref{subsec:Additional-Derivations} shows that 
\begin{equation}
\mathbb{E}_{W}\left[\left\Vert \frac{1}{n}X^{T}W\right\Vert _{\infty}\right]\geq\frac{1}{2}\left(1-\frac{1}{e}\right)\sqrt{\frac{\log L}{4n^{2}}\min_{j,l\in\left\{ 1,...,L\right\} }\sum_{i=1}^{n}\left(X_{ij}-X_{il}\right)^{2}}\label{eq:lead}
\end{equation}
for all $L\geq20$, and 
\begin{equation}
\mathbb{E}_{W}\left[\left\Vert \frac{1}{n}X^{T}W\right\Vert _{\infty}\right]\leq\sqrt{\frac{2\log L}{n^{2}}\max_{j\in\left\{ 1,...,L\right\} }\sum_{i=1}^{n}X_{ij}^{2}}+\sqrt{\frac{8}{n^{2}\log L}\max_{j\in\left\{ 1,...,L\right\} }\sum_{i=1}^{n}X_{ij}^{2}}\label{eq:lead1}
\end{equation}
for all $L\geq2$. While the nonasymptotic validity of our testing
procedure does not require any growth conditions on the dimensionality,
we see from (\ref{eq:lead}) that $\delta_{\alpha,\beta,\infty}$
can tend to zero only when $\frac{\log L}{n}=o(1)$ (if $X$ does
not contain identical columns).

In certain situations, (\ref{eq:3}) and (\ref{eq:ms1}) (with $q=\infty$
in both) can be more conservative than their union bound counterparts.
On the other hand, the extreme example discussed previously suggests
that in situations where $X$ consists of realizations from highly
dependent random vectors, $r_{\alpha,\infty}$ and $\delta_{\alpha,\beta,\infty}$
from the union bound approach can be bigger than (\ref{eq:3}) and
(\ref{eq:ms1}) (with $q=\infty$ in both), respectively, due to the
extra ``$\log L$'' term. This result is similar in spirit to those
of \cite{abr10}, which study bootstrap confidence regions for the
mean of a random vector with the concentration approach. For implementation
of (\ref{eq:p2-1}) with $q=\infty$, we may set $r_{\alpha,\infty}$
to the minimum of (\ref{eq:3}) (with $q=\infty$) and (\ref{eq:bon}).

\subsection{Unknown noise variance }

When no prior information on $\sigma$ is available, $\sqrt{\textrm{Var}\left(Y_{i}\right)}$
may be used as an upper bound. We can easily estimate $\sqrt{\textrm{Var}\left(Y_{i}\right)}$
by $\hat{\sigma}_{Y}=\sqrt{n^{-1}\sum\left(Y_{i}-\bar{Y}\right)^{2}}$
where $\bar{Y}=\frac{1}{n}\sum_{i=1}^{n}Y_{i}$. Proposition 4.1 in
\cite{abr10} implies that 
\[
\sqrt{\textrm{Var}\left(Y_{i}\right)}\leq\left(C_{n}-\frac{1}{\sqrt{n}}\Phi^{-1}\left(\frac{\kappa}{2}\right)\right)^{-1}\hat{\sigma}_{Y}:=\bar{B}_{\kappa}
\]
with probability at least $1-\kappa$, where $C_{n}=\sqrt{\frac{2}{n}}\frac{\Gamma\left(n/2\right)}{\Gamma\left((n-1)/2\right)}=1-O(n^{-1})$. 

In problems where $\textrm{Var}(W_{i})$ is a constant over $i$,
$V$ is fixed, and the only source of randomness in $Y$ comes from
$W$, replacing $\sigma$ with $\bar{B}_{\kappa}$ does not make $r_{\alpha,q}$
a more conservative threshold for constructing confidence regions.
In problems with a random design, using $\bar{B}_{\kappa}$ could
result in confidence regions that are more conservative. We note that
it is rather challenging to estimate $\sigma$ precisely and obtain
a sharp threshold simultaneously within our framework. The main reason
is that our guarantee does not require a small $\sqrt{\frac{1}{n}\sum_{i=1}^{n}\left[g\left(V_{i};\theta^{*}\right)-g\left(V_{i};\hat{\theta}\right)\right]^{2}}$
with high probability, which seems to be needed for consistent estimation
of $\sigma$. On the other hand, if we were able to ensure a small
prediction error, more assumptions might be required and our nonasymptotic
control is likely to involve unknown nuisance parameters that are
hard to estimate.

\section{Simulations}

In this section, we illustrate the performance of our procedure through
a Monte-Carlo experiment based on $100$ repetitions. The data are
generated according to the following model: 
\[
Y_{i}=\sum_{l=1}^{k}v_{il}\alpha_{l}^{*}+\gamma^{*}\exp\left(\sum_{l=1}^{k}v_{il}\tau_{l}^{*}\right)+W_{i},\qquad i=1,...,n.
\]
The matrix $v$ consists of $n$ rows, which are fixed realizations
of i.i.d. draws from the normal distribution $\mathcal{N}\left(\mathbf{0}_{k},\,\Sigma\right)$
where $\Sigma_{jj}=1$ and $\Sigma_{jj^{'}}=0.5$ for $j\neq j^{'}$,
$j,j^{'}\in\left\{ 1,...,k\right\} $. The $\left(i,\,l\right)$th
entry of $v$ is denoted by $v_{il}$; $\alpha^{*}=\left\{ \alpha_{l}^{*}\right\} _{l=1}^{k}\in\mathbb{R}^{k}$,
$\tau^{*}=\left\{ \tau_{l}^{*}\right\} _{l=1}^{k}\in\mathbb{R}^{k}$,
and $\gamma^{*}\in\mathbb{R}$ are the unknown coefficients (as a
result, $p=2k+1$). For each of the $100$ repetitions, the noise
vector $W$ is randomly drawn from $\mathcal{N}(\mathbf{0}_{n},\,\sigma^{2}\boldsymbol{I}_{n})$
with $\sigma=0.5$.

Let $v_{i}$ denote the $i$th row of $v$ and $g\left(v_{i};\theta^{*}\right)=\sum_{l=1}^{k}v_{il}\alpha_{l}^{*}+\gamma^{*}\exp\left(\sum_{l=1}^{k}v_{il}\tau_{l}^{*}\right)$.
Our hypotheses take the form 
\begin{align}
H_{0}:\,\frac{1}{n}\sum_{i=1}^{n}\frac{\partial g\left(v_{i};\theta^{*}\right)}{\partial v_{il}} & \in\left[0,\,0.8\right],\qquad\forall l\in M\subseteq\left\{ 1,...,k\right\} \label{eq:h2-1}\\
H_{1}:\,\frac{1}{n}\sum_{i=1}^{n}\frac{\partial g\left(v_{i};\theta^{*}\right)}{\partial v_{il}} & \notin\left[0,\,0.8\right],\qquad\forall l\in M\subseteq\left\{ 1,...,k\right\} \nonumber 
\end{align}
where 
\[
\frac{1}{n}\sum_{i=1}^{n}\frac{\partial g\left(v_{i};\theta^{*}\right)}{\partial v_{il}}=\alpha_{l}^{*}+\gamma^{*}\tau_{l}^{*}\frac{1}{n}\sum_{i=1}^{n}\exp\left(\sum_{l=1}^{k}v_{il}\tau_{l}^{*}\right)
\]
is the APE of $v_{il}$s on $\mathbb{E}\left(Y_{i}\right)=g\left(v_{i};\theta^{*}\right)$.
Two samples sizes, $n=\left\{ 30,\,90\right\} $, are considered in
our experiment. We first look at a low-dimensional scenario where
$k=1$ (i.e., $v_{i}=v_{i1}$) and $\theta^{*}$ consists of $3$
unknown components ($p=3$). For this scenario, we have $M=\left\{ 1\right\} $
in (\ref{eq:h2-1}) ($m=1$). We then look at a high-dimensional scenario
where $k=15$ (i.e., $v_{i}=\left\{ v_{il}\right\} _{l=1}^{15}$)
and $\theta^{*}$ consists of $31$ unknown components ($p=31$).
For this scenario, we consider $M=\left\{ 1,...,k\right\} $ in (\ref{eq:h2-1})
($m=15$). To examine the impact of $L$, we look at the choice $X_{i}=\left(v_{i1},\,v_{i1}^{2},\,v_{i1}^{3}\right)\in\mathbb{R}^{3}$
versus the choice $X_{i}=\left(v_{i1},\,v_{i1}^{2},\,v_{i1}^{3},\,v_{i1}^{4}\right)\in\mathbb{R}^{4}$
in the low-dimensional scenario, and $X_{i}=\left\{ v_{ij},\,v_{ij}^{2},\,v_{ij}^{3}\right\} _{j=1}^{15}\in\mathbb{R}^{45}$
versus $X_{i}=\left\{ v_{ij},\,v_{ij}^{2},\,v_{ij}^{3},\,v_{ij}^{4}\right\} _{j=1}^{15}\in\mathbb{R}^{60}$
in the high-dimensional scenario. The entries in $X$ are normalized
such that $\frac{1}{n}\sum_{i=1}^{n}X_{ij}^{2}=1$ for all $j=1,...,L$.

We apply program (\ref{eq:p2-1}) with $q=\infty$ and set $r_{\alpha,\infty}$
to the minimum of (\ref{eq:3}) (with $q=\infty$) and (\ref{eq:bon}).
For each of the $100$ repetitions, we take $R=10000$ i.i.d. draws
($Z_{r}$s) from $\mathcal{N}\left(\mathbf{0}_{n},\,\boldsymbol{I}_{n}\right)$
and choose $\alpha_{1}=0.049$, $\alpha_{2}=0.001$ (i.e., $\alpha=0.05$)
to balance between $\tau_{\alpha_{1},\infty}$ and $\sqrt{\frac{1}{R}}\tau_{\alpha_{2},\infty}$
in (\ref{eq:3}). For (\ref{eq:bon}), we simply choose $\alpha=0.05$.
We compare the actual separation $\left\Vert \frac{1}{n}\sum_{i=1}^{n}X_{i}\left[g\left(V_{i};\theta^{*}\right)-g\left(V_{i};\hat{\theta}_{\alpha}\right)\right]\right\Vert _{\infty}$
with $2r_{\alpha,\infty}$. By setting $\beta_{1}=0.001$, $\beta_{2}=0.049$
in (\ref{eq:ms1}) (and $\beta=0.05$ in (\ref{eq:ms-bon})), note
that $2r_{\alpha,\infty}$ approximates (respectively, coincides with)
$\delta_{\alpha,\beta,\infty}$. Tables 1 and 2 exhibit:
\begin{lyxlist}{00.00.0000}
\item [{\textbf{(i)}}] $\frac{1}{n}\sum_{i=1}^{n}\frac{\partial g\left(v_{i};\theta^{*}\right)}{\partial v_{il}}$
for $l\in M$, 
\item [{\textbf{(ii)}}] the average of $\left\Vert \frac{1}{n}\sum_{i=1}^{n}X_{i}\left[g\left(V_{i};\theta^{*}\right)-g\left(V_{i};\hat{\theta}_{\alpha}\right)\right]\right\Vert _{\infty}$
over $100$ repetitions, 
\item [{\textbf{(iii)}}] the average of $2r_{\alpha,\infty}$ over $100$
repetitions, 
\item [{\textbf{(iv)}}] the coverage probability, 
\item [{\textbf{(v)}}] the rejection probability. 
\end{lyxlist}
The evidence from our simulation study supports the main points of
Section \ref{sec:Gaussian-Regressions}. For our procedure to reject
the null hypothesis, all it takes is sufficient separation in $\left\Vert \frac{1}{n}\sum_{i=1}^{n}X_{i}\left[g\left(V_{i};\theta^{*}\right)-g\left(V_{i};\hat{\theta}_{\alpha}\right)\right]\right\Vert _{\infty}$
and $\theta^{*}$ need not have any form of sparsity. Also in view
of (\ref{eq:7b}) and (\ref{eq:7b-low}), it is not surprising that
the coverage probabilities of our procedure are not affected by whether
$\theta^{*}$ is sparse or not. 

In contrast to ``undercoverage'' commonly reported in many asymptotic
procedures, the coverage probabilities shown in Tables 1 and 2 suggest
that our method tends to be conservative. The actual separation (ii)
needed to achieve a power of around $95\%$ or higher is somewhat
comparable to (and smaller than) the theoretical prediction (iii).
This result is plausible given that our control on the Type II error
only states that $\beta$ is an upper bound on the probability of
our procedure failing to reject $H_{0}$. 

As shown in Tables 1 and 2, the values of $\frac{1}{n}\sum_{i=1}^{n}\frac{\partial g\left(v_{i};\theta^{*}\right)}{\partial v_{il}}$s
(i) that make most of the rejection probabilities around $95\%$ or
higher decrease as $n$ increases from $30$ to $90$. This finding
is intuitive: keeping all the other factors the same, the numbers
of constraints and ``free'' parameters to be determined in (\ref{eq:p2-1})
remain unchanged but $r_{\alpha,\infty}$ decreases as $n$ increases.
Consequently, it takes smaller $\frac{1}{n}\sum_{i=1}^{n}\frac{\partial g\left(v_{i};\theta^{*}\right)}{\partial v_{il}}$s
for our procedure to correctly reject $H_{0}$. In addition, a larger
$L$ improves the power in most cases, for which we have provided
an explanation in Section \ref{subsec:Discussions}.\\
\\

\begin{center}
\begin{tabular}{|c|c|c|c|c|c|}
\hline 
\multicolumn{6}{|c|}{\textbf{\small{}Table 1:}{\small{} $p=3$, $M=\left\{ 1\right\} $}}\tabularnewline
\hline 
 & \textbf{\small{}i} & \textbf{\small{}ii} & \textbf{\small{}iii} & \textbf{\small{}iv} & \textbf{\small{}v}\tabularnewline
\hline 
\textbf{\small{}$n=30$, $L=4$} & \textbf{\small{}$1.482$} & \textbf{\small{}$0.091$} & \textbf{\small{}$0.106$} & \textbf{\small{}$1$} & \textbf{\small{}$93\%$}\tabularnewline
\hline 
\textbf{\small{}$n=30$, $L=3$} & \textbf{\small{}$1.482$} & \textbf{\small{}$0.068$} & \textbf{\small{}$0.103$} & \textbf{\small{}$1$} & \textbf{\small{}$64\%$}\tabularnewline
\hline 
\textbf{\small{}$n=90$, $L=4$} & \textbf{\small{}$1.159$} & \textbf{\small{}$0.028$} & \textbf{\small{}$0.035$} & \textbf{\small{}$1$} & \textbf{\small{}$93\%$}\tabularnewline
\hline 
\textbf{\small{}$n=90$, $L=3$} & \textbf{\small{}$1.159$} & \textbf{\small{}$0.028$} & \textbf{\small{}$0.034$} & \textbf{\small{}$1$} & \textbf{\small{}$96\%$}\tabularnewline
\hline 
\multicolumn{6}{c}{{\footnotesize{}For $n=30$, $\alpha^{*}=\tau^{*}=0.657$, $\gamma^{*}=1$}}\tabularnewline
\multicolumn{6}{c}{{\footnotesize{}For $n=90$, $\alpha^{*}=\tau^{*}=0.533$, $\gamma^{*}=1$}}\tabularnewline
\end{tabular}\\
\medskip{}
\medskip{}
\par\end{center}

\begin{center}
\begin{tabular}{|c|c|c|c|c|c|}
\hline 
\multicolumn{6}{|c|}{\textbf{\small{}Table 2: }{\small{}$p=31$, $M=\left\{ 1,...,15\right\} $}}\tabularnewline
\hline 
 & \textbf{\small{}i} & \textbf{\small{}ii} & \textbf{\small{}iii} & \textbf{\small{}iv} & \textbf{\small{}v}\tabularnewline
\hline 
\textbf{\small{}$n=30$, $L=60$} & {\small{}$\left(1.735,...,1.735\right)$} & {\small{}$0.099$} & {\small{}$0.132$} & {\small{}$1$} & {\small{}$100\%$}\tabularnewline
\hline 
\textbf{\small{}$n=30$, $L=45$} & {\small{}$\left(1.735,...,1.735\right)$} & {\small{}$0.096$} & {\small{}$0.129$} & {\small{}$1$} & {\small{}$95\%$}\tabularnewline
\hline 
\textbf{\small{}$n=90$, $L=60$} & {\small{}$\left(0.929,...,0.929\right)$} & {\small{}$0.035$} & {\small{}$0.044$} & {\small{}$1$} & {\small{}$99\%$}\tabularnewline
\hline 
\textbf{\small{}$n=90$, $L=45$} & {\small{}$\left(0.929,...,0.929\right)$} & {\small{}$0.033$} & {\small{}$0.043$} & {\small{}$1$} & {\small{}$94\%$}\tabularnewline
\hline 
\multicolumn{6}{c}{{\footnotesize{}For $n=30$, $\alpha^{*}=\tau^{*}=\left(0.194,...,0.194\right)$,
$\gamma^{*}=1$}}\tabularnewline
\multicolumn{6}{c}{{\footnotesize{}For $n=90$, $\alpha^{*}=\tau^{*}=\left(0.172,...,0.172\right)$,
$\gamma^{*}=1$}}\tabularnewline
\end{tabular}\\
\medskip{}
\par\end{center}

\begin{center}
\medskip{}
\par\end{center}

\section{Extensions }

Beyond Gaussian regressions, it is possible to establish some nonasymptotic
justifications for inference in statistical models that involve non-Gaussian
responses. 

\subsection{Regressions with non-Gaussian noise}

Our analyses in Section \ref{sec:Gaussian-Regressions} exploit sharp
concentration of Lipschitz functions of Gaussian variables. These
analyses can be extended to regression models where the noise vector
$W$ is either bounded or has a strongly log-concave distribution.
In particular, we have the following analogues of (\ref{eq:prop2-1}).\\
\textbf{\textit{}}\\
\textbf{\textit{Lemma 4.1}}\textit{. Suppose $W$ has a strongly log-concave}\footnote{A strongly log-concave distribution is a distribution with density
$\mathsf{p}(z)=\exp\left(-\psi(z)\right)$ such that for some $\varphi>0$
and all $\lambda\in[0,\,1]$, $z,\,z^{'}\in\mathbb{R}^{n}$, $\lambda\psi(z)+(1-\lambda)\psi(z^{'})-\psi(\lambda z+(1-\lambda)z^{'})\geq\frac{\varphi}{2}\lambda(1-\lambda)\left\Vert z-z^{'}\right\Vert _{2}^{2}.$}\textit{ distribution with parameter $\varphi$. Then for any $q\in\left[1,\,\infty\right]$,
we have 
\begin{equation}
\mathbb{P}\left\{ \left\Vert \frac{1}{n}X^{T}W\right\Vert _{q}\geq\mathbb{E}\left[\left\Vert \frac{1}{n}X^{T}W\right\Vert _{q}\right]+t\right\} \leq\exp\left(\frac{-n\varphi t^{2}}{2\left\Vert \sqrt{\frac{1}{n}\sum_{i=1}^{n}X_{i}^{2}}\right\Vert _{q}^{2}}\right).\label{eq:lemma4.1}
\end{equation}
}\textbf{Remarks}. Let $g\left(V;\theta^{*}\right):=\left\{ g\left(V_{i};\theta^{*}\right)\right\} _{i=1}^{n}$.
For a fixed design $V$, if $Y\sim\mathcal{N}\left(g\left(V;\theta^{*}\right),\,\Sigma\right)$
and $\Sigma\succ0$, $\varphi$ can be set to the smallest eigenvalue
of $\Sigma^{-1}$. Beyond a normal distribution, \cite{sw14} discuss
quite a few examples of strongly log-concave distributions.\textbf{\textit{
}}\\
\textbf{\textit{}}\\
\textbf{\textit{Lemma 4.2}}\textit{. Suppose $W$ consists of independent
random variables, all of which are supported on $\left[a,\,b\right]$.
Then for any $q\in\left[1,\,\infty\right]$, we have 
\begin{equation}
\mathbb{P}\left\{ \left\Vert \frac{1}{n}X^{T}W\right\Vert _{q}\geq\mathbb{E}\left[\left\Vert \frac{1}{n}X^{T}W\right\Vert _{q}\right]+t\right\} \leq\exp\left(\frac{-nt^{2}}{2\left(b-a\right)^{2}\left\Vert \sqrt{\frac{1}{n}\sum_{i=1}^{n}X_{i}^{2}}\right\Vert _{q}^{2}}\right).\label{eq:lemma4.2}
\end{equation}
}If we know the distribution of $W$, our analyses from Section \ref{sec:Gaussian-Regressions}
can be, in principle, extended to construct inference procedures for
regression models where $W$ is either bounded or has a strongly log-concave
distribution. However, sometimes we might not know the distribution
for $W$; instead, we may have more information on the distribution
of $Y$ than the distribution of $W$. In some applications, we might
know $Y$ consists of entries supported on $\left[a,\,b\right]$.
For example, \cite{pw08} and \cite{wz17} estimate the APE of spending
on math pass rates ($Y_{i}\in\left[0,\,1\right]$) under the assumption
$\mathbb{E}\left(Y_{i}\vert V_{i}\right)=\Phi\left(V_{i}\theta^{*}\right)$,
where $\Phi(\cdot)$ denotes the standard normal c.d.f. and $V_{i}$
include the spending variable as well as other covariates. Another
example is the binary response model
\begin{equation}
\mathbb{P}\left(Y_{i}=1\vert V_{i}\right)=g\left(V_{i};\theta^{*}\right),\qquad i=1,...,n,\label{eq:binary}
\end{equation}
where $Y_{i}\in\left\{ 0,\,1\right\} $ and the functional form of
$g\left(V_{i};\theta^{*}\right)$ is known; for example, $g$ may
be a ``probit'' or a ``logit'' in (\ref{eq:binary}) and $g\left(V_{i};\theta^{*}\right)=g\left(V_{i}\theta^{*}\right)$.
Under the assumption 
\begin{equation}
\mathbb{E}\left(Y_{i}\vert V_{i}\right)=g\left(V_{i};\theta^{*}\right),\label{eq:model}
\end{equation}
both binary and bounded response models can be treated in the same
framework. 

\subsection{Bounded responses }

In what follows, we consider (\ref{eq:model}) where $a\leq Y_{i}\leq b$
for all $i$, the functional form of $g\left(V_{i};\theta^{*}\right)$
is known and possibly nonlinear in $\theta^{*}$. Without loss of
generality, we assume $a=0$ and $b=1$. Throughout this section,
we use $\mathbb{E}_{Y|V}\left[\cdot\right]$ to denote the expectation
over the distribution of $Y$ conditioning on $V$; for an i.i.d.
sequence of Radamacher random variables, $\varepsilon=\left\{ \varepsilon_{i}\right\} _{i=1}^{n}$
(independent of $Y$ and $V$), we use $\mathbb{E}_{\varepsilon}\left[\cdot\right]$
to denote the expectation over $\varepsilon$ only, conditioning on
$Y$ and $V$, and $\mathbb{E}_{\varepsilon,Y|V}\left[\cdot\right]$
to denote the expectation over the distribution of $(\varepsilon,\,Y)$
conditioning on $V$.

As in Section \ref{subsec:The-Bonferroni-alternative}, for $q=\infty$,
the choice of $r_{\alpha,\infty}$ can be based on simple union bounds.
We omit discussions of this strategy but focus on the concentration
approach for the more general cases $q\in[1,\,\infty]$. Like in the
regression problem, we first establish the concentration of 
\[
\left\Vert \frac{1}{n}\sum_{i=1}^{n}X_{i}\left[Y_{i}-g\left(V_{i};\theta^{*}\right)\right]\right\Vert _{q}
\]
around its expectation 
\begin{equation}
S_{\theta^{*}}:=\mathbb{E}_{Y|V}\left[\left\Vert \frac{1}{n}\sum_{i=1}^{n}X_{i}\left[Y_{i}-g\left(V_{i};\theta^{*}\right)\right]\right\Vert _{q}\right].\label{eq:exp}
\end{equation}
Previously we have simply replaced $\mathbb{E}_{W}\left[\left\Vert \frac{1}{n}X^{T}W\right\Vert _{q}\right]$
in (\ref{eq:prop2-1}) with its Monte Carlo approximation $\frac{\sigma}{R}\sum_{r=1}^{R}\left\Vert \frac{1}{n}X^{T}Z_{r}\right\Vert _{q}$
and a ``small'' deviation. This strategy cannot be applied to the
expectation $S_{\theta^{*}}$ directly. Instead, we first seek a reasonable
upper bound which involves only $\left\{ Y,\,X\right\} $ and random
variables from a known distribution. This result is stated in the
following proposition. \\
 \\
\textbf{\textit{Proposition 4.1}}\textit{. Assume $Y=\left\{ Y_{i}\right\} _{i=1}^{n}$
consists of independent random variables, where $0\leq Y_{i}\leq1$
for all $i$ (w.l.o.g.). For any $q\in\left[1,\,\infty\right]$, we
have 
\begin{equation}
\mathbb{P}\left\{ \left\Vert \frac{1}{n}\sum_{i=1}^{n}X_{i}\left[Y_{i}-g\left(V_{i};\theta^{*}\right)\right]\right\Vert _{q}\geq S_{\theta^{*}}+t\right\} \leq\exp\left(\frac{-nt^{2}}{2\left\Vert \sqrt{\frac{1}{n}\sum_{i=1}^{n}X_{i}^{2}}\right\Vert _{q}^{2}}\right).\label{eq:prop3.1}
\end{equation}
Let $\varepsilon=\left\{ \varepsilon_{i}\right\} _{i=1}^{n}$ be an
i.i.d. sequence of Radamacher random variables independent of $Y$
and $V$. Under (\ref{eq:model}), we have 
\begin{equation}
\mathbb{E}_{\varepsilon,Y|V}\left\{ \left\Vert \frac{1}{2n}\sum_{i=1}^{n}\varepsilon_{i}X_{i}\left[Y_{i}-g\left(V_{i};\theta^{*}\right)\right]\right\Vert _{q}\right\} \leq S_{\theta^{*}}\leq2\mathbb{E}_{\varepsilon,Y|V}\left\{ \left\Vert \frac{1}{n}\sum_{i=1}^{n}\varepsilon_{i}X_{i}Y_{i}\right\Vert _{q}\right\} .\label{eq:sandwich}
\end{equation}
}\textbf{Remarks}. Note that bound (\ref{eq:prop3.1}) holds for any
fixed $\theta$ (not just the true coefficient vector, $\theta^{*}$).
However, (\ref{eq:sandwich}) relies crucially on the model assumption
(\ref{eq:model}).

The upper bound in (\ref{eq:sandwich}) can be viewed as the symmetrized
version of $S_{\theta^{*}}$. Considering a collection of i.i.d. Radamacher
random draws (independent of $Y$ and $V$), \textit{ 
\begin{equation}
\left\{ \varepsilon_{ir}:\,i=1,...,n,\,r=1,...,R\right\} ,\label{eq:rad}
\end{equation}
}we can replace $S_{\theta^{*}}$ with $\frac{2}{R}\sum_{r=1}^{R}\left\Vert \frac{1}{n}\sum_{i=1}^{n}\varepsilon_{ir}Y_{i}X_{i}\right\Vert _{q}$
(a Monte-Carlo approximation of the symmetrized version) and some
``small'' deviations. The complementary lower bound in (\ref{eq:sandwich})
suggests that $S_{\theta^{*}}$ and its symmetrized version have the
same magnitude. 

The first step is to relate $S_{\theta^{*}}$ with $\frac{2}{R}\sum_{r=1}^{R}\left\Vert \frac{1}{n}\sum_{i=1}^{n}\varepsilon_{ir}Y_{i}X_{i}\right\Vert _{q}$
as shown in the following proposition. \\
\\
\textbf{\textit{Proposition 4.2}}\textit{. Assume (\ref{eq:model})
where $Y=\left\{ Y_{i}\right\} _{i=1}^{n}$ consists of independent
random variables, $0\leq Y_{i}\leq1$ for all $i$ (w.l.o.g.). Given
(\ref{eq:rad}) which is independent of $Y$ and $V$, for any $q\in\left[1,\,\infty\right]$,
we have 
\begin{equation}
\left\Vert \frac{1}{n}\sum_{i=1}^{n}X_{i}\left[Y_{i}-g\left(V_{i};\theta^{*}\right)\right]\right\Vert _{q}\geq\frac{2}{R}\sum_{r=1}^{R}\left\Vert \frac{1}{n}\sum_{i=1}^{n}\varepsilon_{ir}Y_{i}X_{i}\right\Vert _{q}+t_{1}+2t_{2}+2t_{3}\label{eq:c1-2}
\end{equation}
with probability no greater than $\alpha\in(0,\,1)$, where  
\begin{eqnarray*}
t_{1} & = & \tau_{\alpha_{1},q}=\left\Vert \sqrt{\frac{1}{n}\sum_{i=1}^{n}X_{i}^{2}}\right\Vert _{q}\sqrt{\frac{2}{n}\log\frac{1}{\alpha_{1}}},\\
t_{2} & = & \tau_{\alpha_{2},q}=\left\Vert \sqrt{\frac{1}{n}\sum_{i=1}^{n}X_{i}^{2}}\right\Vert _{q}\sqrt{\frac{2}{n}\log\frac{1}{\alpha_{2}}},\\
t_{3} & = & \frac{2}{\sqrt{R}}\tau_{\alpha_{3},q}=\left\Vert \sqrt{\frac{1}{n}\sum_{i=1}^{n}X_{i}^{2}}\right\Vert _{q}\sqrt{\frac{8}{nR}\log\frac{1}{\alpha_{3}}},
\end{eqnarray*}
for chosen $\alpha_{1},\,\alpha_{2},\,\alpha_{3}>0$ such that $\sum_{k=1}^{3}\alpha_{k}=\alpha$.}\\

Based on (\ref{eq:c1-2}) along with the choices of $t_{1}$, $t_{2}$
and $t_{3}$ above, we set in (\ref{eq:p2-1}), 
\begin{equation}
r_{\alpha,q}=\frac{2}{R}\sum_{r=1}^{R}\left\Vert \frac{1}{n}\sum_{i=1}^{n}\varepsilon_{ir}Y_{i}X_{i}\right\Vert _{q}+\tau_{\alpha_{1},q}+2\tau_{\alpha_{2},q}+\frac{4}{\sqrt{R}}\tau_{\alpha_{3},q}.\label{eq:3-1}
\end{equation}
Under $H_{0}$, $\left(\theta^{*},\,0\right)$ is an optimal solution
to (\ref{eq:p2-1}) with $r_{\alpha,q}$ specified in (\ref{eq:3-1}).
Consequently, a (practical) optimal solution to (\ref{eq:p2-1}) must
satisfy \textit{ 
\begin{equation}
\mathbb{P}_{0}\left\{ \Psi_{q}(\hat{\theta}_{\alpha})\geq r_{\alpha,q}\right\} \leq\alpha.\label{eq:13-3}
\end{equation}
}Moreover, following the argument used for deriving (\ref{eq:low2})
and (\ref{eq:6-1}), we can show that, an optimal solution $\left(\hat{\theta}_{\alpha},\,\hat{\mu}_{\alpha}\right)$
to (\ref{eq:p2-1}) must satisfy 
\begin{eqnarray}
\left\Vert \frac{1}{n}\sum_{i=1}^{n}X_{i}\left[g\left(V_{i};\theta^{*}\right)-g\left(V_{i};\hat{\theta}_{\alpha}\right)\right]\right\Vert _{q} & \geq & \hat{\mu}_{\alpha},\label{eq:CI_11}\\
\left\Vert \frac{1}{n}\sum_{i=1}^{n}X_{i}\left[g\left(V_{i};\theta^{*}\right)-g\left(V_{i};\hat{\theta}_{\alpha}\right)\right]\right\Vert _{q} & \leq & 2r_{\alpha,q}+\hat{\mu}_{\alpha},\label{eq:CI_22}
\end{eqnarray}
with probability at least $1-\alpha$. 

For general $q\in\left[1,\,\infty\right]$, the strategy where we
replace $S_{\theta^{*}}$ in (\ref{eq:exp}) by 
\[
\frac{2}{R}\sum_{r=1}^{R}\left\Vert \frac{1}{n}\sum_{i=1}^{n}\varepsilon_{ir}Y_{i}X_{i}\right\Vert _{q}
\]
plus some deviations only requires the correct specification of the
conditional mean of $Y_{i}$, i.e., (\ref{eq:model}). This treatment
delivers generic confidence regions in the form of (\ref{eq:CI_11})
and (\ref{eq:CI_22}) for binary and fractional responses. 

\subsection*{Acknowledgment}

I am grateful to the associate editor's insightful comments. I also
thank Guang Cheng at Purdue University for discussions that improved
the clarity of an earlier version (arXiv:1808.07127v1) and pointing
out several references. All errors are my own.

\appendix

\section{Supplementary materials}

\subsection{Alternative formulation \label{subsec:Alternative-formulation}}

Instead of (\ref{eq:p2-1}), we can work with an alternative formulation:
\begin{align}
\left(\hat{\theta}_{\alpha},\,\hat{\mu}_{\alpha}\right)\in\arg\min_{(\theta_{\alpha},\mu_{\alpha})\in\mathbb{R}^{p}\times\mathbb{R}^{p}}\left\Vert \mu_{\alpha}\right\Vert _{\tilde{q}}\nonumber \\
\textrm{subject to: }\left\Vert \frac{1}{n}\sum_{i=1}^{n}X_{i}\left[Y_{i}-g\left(V_{i};\theta_{\alpha}\right)\right]-\mu_{\alpha}\right\Vert _{q}\leq r_{\alpha,q},\label{eq:p1-1}\\
h(\theta_{\alpha})\in\Omega.\nonumber 
\end{align}
Here we slightly abuse the notations, where $\hat{\mu}_{\alpha}$
(also $\mu_{\alpha}$) in (\ref{eq:p1-1}) is a vector and in (\ref{eq:p2-1})
is a scalar. Like in (\ref{eq:p2-1}), we suppress the dependence
of $(\hat{\theta}_{\alpha},\hat{\mu}_{\alpha})$ in (\ref{eq:p1-1})
on $\left(q,\,\tilde{q}\right)$ for notational simplicity. 

To compare (\ref{eq:p1-1}) with (\ref{eq:p2-1}) from the computational
perspective, we let $\mathcal{F}_{1}^{\alpha}$ denote the set of
$(\theta_{\alpha},\mu_{\alpha})$ that are feasible for (\ref{eq:p1-1})
and $\mathcal{F}_{1,\theta}^{\alpha}$ denote the set of $\theta_{\alpha}$
from $\mathcal{F}_{1}^{\alpha}$; similarly, $\mathcal{F}_{2}^{\alpha}$
and $\mathcal{F}_{2,\theta}^{\alpha}$ are defined with regard to
(\ref{eq:p2-1}). Note that an element $(\tilde{\theta}_{\alpha},\tilde{\mu}_{\alpha})$
in $\mathcal{F}_{1}^{\alpha}$ implies 
\[
\left\Vert \frac{1}{n}\sum_{i=1}^{n}X_{i}\left[Y_{i}-g\left(V_{i};\tilde{\theta}_{\alpha}\right)\right]\right\Vert _{q}\leq r_{\alpha,q}+\left\Vert \tilde{\mu}_{\alpha}\right\Vert _{q};
\]
that is, $(\tilde{\theta}_{\alpha},\left\Vert \tilde{\mu}_{\alpha}\right\Vert _{q})\in\mathcal{F}_{2}^{\alpha}$.
Consequently, $\mathcal{F}_{1,\theta}^{\alpha}\subseteq\mathcal{F}_{2,\theta}^{\alpha}$.
On the other hand, the objective function in (\ref{eq:p1-1}) is minimized
over an $L-$dimensional vector as opposed to a scalar in (\ref{eq:p2-1}).
These facts suggest that the choice between (\ref{eq:p1-1}) and (\ref{eq:p2-1})
incurs some trade-offs in terms of computational cost. 

Theorem 2.1 holds for the alternative formulation. Moreover, the following
results exhibit the ``ideal'' and practical confidence regions.

\subsubsection{Ideal confidence regions }

\textbf{\textit{Proposition A.2.1}}\textit{. Assume (\ref{eq:reg})
where $W\sim\mathcal{N}(\mathbf{0}_{n},\,\sigma^{2}\boldsymbol{I}_{n})$
and is independent of $V$. Then, an optimal solution $\left(\hat{\theta}_{\alpha}^{*},\,\hat{\mu}_{\alpha}^{*}\right)$
to (\ref{eq:p1-1}) must satisfy 
\begin{eqnarray}
\left\Vert \frac{1}{n}\sum_{i=1}^{n}X_{i}\left[g\left(V_{i};\theta^{*}\right)-g\left(V_{i};\hat{\theta}_{\alpha}^{*}\right)\right]\right\Vert _{\tilde{q}} & \geq & \left\Vert \hat{\mu}_{\alpha}^{*}\right\Vert _{\tilde{q}},\label{eq:low1-1}\\
\left\Vert \frac{1}{n}\sum_{i=1}^{n}X_{i}\left[g\left(V_{i};\theta^{*}\right)-g\left(V_{i};\hat{\theta}_{\alpha}^{*}\right)\right]-\hat{\mu}_{\alpha}^{*}\right\Vert _{q} & \leq & 2r_{\alpha,q}^{*},\label{eq:6-2}
\end{eqnarray}
with probability at least $1-\alpha$, where $r_{\alpha,q}^{*}$ is
specified in (\ref{eq:8}). }

\subsubsection{Practical confidence regions }

Let $\left(\hat{\theta}_{\alpha},\,\hat{\mu}_{\alpha}\right)$ be
an optimal solution to (\ref{eq:p1-1}) with $r_{\alpha,q}$ specified
in (\ref{eq:3}). Then we have

\begin{eqnarray}
 &  & \left\Vert \frac{1}{n}\sum_{i=1}^{n}X_{i}\left[g\left(V_{i};\theta^{*}\right)-g\left(V_{i};\hat{\theta}_{\alpha}\right)\right]-\hat{\mu}_{\alpha}\right\Vert _{q}\nonumber \\
 & \leq & \left\Vert \frac{1}{n}\sum_{i=1}^{n}X_{i}\left[Y_{i}-g\left(V_{i};\hat{\theta}_{\alpha}\right)\right]-\hat{\mu}_{\alpha}\right\Vert _{q}+\left\Vert \frac{1}{n}X^{T}W\right\Vert _{q}\nonumber \\
 & \leq & \frac{2\sigma}{R}\sum_{r=1}^{R}\left\Vert \frac{1}{n}X^{T}Z_{r}\right\Vert _{q}+2\tau_{\alpha_{1},q}+2\sqrt{\frac{1}{R}}\tau_{\alpha_{2},q}\label{eq:7}
\end{eqnarray}
and 
\begin{equation}
\left\Vert \frac{1}{n}\sum_{i=1}^{n}X_{i}\left[g\left(V_{i};\theta^{*}\right)-g\left(V_{i};\hat{\theta}_{\alpha}\right)\right]\right\Vert _{\tilde{q}}\geq\left\Vert \hat{\mu}_{\alpha}\right\Vert _{\tilde{q}}\label{eq:7-low}
\end{equation}
with probability at least $1-\alpha$. 

\subsection{Preliminary \label{subsec:Preliminary}}

Here we include several classical results which are used in the main
proofs. We first introduce a definition of sub-Gaussian variables.
\\
 \\
\textbf{\textit{Definition A.1}}\textit{. A zero-mean random variable
$U_{1}$ is sub-Gaussian if there is a $\nu>0$ such that 
\begin{equation}
\mathbb{E}\left[\exp\left(\lambda U_{1}\right)\right]\leq\exp\left(\frac{\lambda^{2}\nu^{2}}{2}\right)\label{eq:sub}
\end{equation}
for all $\lambda\in\mathbb{R}$, and we refer to $\nu$ as the sub-Gaussian
parameter. }\\
 \\
\textbf{Remarks}. 
\begin{enumerate}
\item Using the Chernoff bound, one can show that any zero-mean random variable
$U_{1}$ obeying (\ref{eq:sub}) satisfies 
\begin{eqnarray}
\mathbb{P}\left(U_{1}\leq-t\right) & \leq & \exp\left(-\frac{t^{2}}{2\nu^{2}}\right),\label{eq:lower}\\
\mathbb{P}\left(U_{1}\geq t\right) & \leq & \exp\left(-\frac{t^{2}}{2\nu^{2}}\right),\label{eq:upper}
\end{eqnarray}
for all $t\geq0$. 
\item Let $\left\{ U_{i}\right\} _{i=1}^{R}$ be independent zero-mean sub-Gaussian
random variables, each with parameter at most $\nu$. Then $R^{-1}\sum_{i=1}^{R}U_{i}$
is sub-Gaussian with parameter at $\nu/\sqrt{R}$. To see this, note
that for all $\lambda\in\mathbb{R}$,
\begin{eqnarray}
\mathbb{E}\left[\exp\left(\frac{\lambda}{R}\sum_{i=1}^{R}U_{i}\right)\right] & = & \prod_{i=1}^{R}\mathbb{E}\left[\exp\left(\frac{\lambda U_{i}}{R}\right)\right]\nonumber \\
 & \leq & \prod_{i=1}^{R}\exp\left(\frac{\lambda^{2}\nu^{2}}{2R^{2}}\right)\nonumber \\
 & = & \exp\left(\frac{\lambda^{2}\nu^{2}}{2R}\right).\label{eq:subn}
\end{eqnarray}
\end{enumerate}
The following result exhibits the type of sub-Gaussian variables that
are of interest to our analysis. \\
 \\
\textbf{\textit{Lemma A.1}}\textit{.} \textit{Suppose $U=\left\{ U_{i}\right\} _{i=1}^{n}$
has a strongly log-concave distribution with parameter $\varphi>0$
and $f:\,\mathbb{R}^{n}\rightarrow\mathbb{R}$ is $\mathcal{L}-$Lipschitz
with respect to the Euclidean norm. Then for all $\lambda\in\mathbb{R}$,
we have
\begin{equation}
\mathbb{E}\left[\exp\left(\lambda\left\{ f(U)-\mathbb{E}\left[f(U)\right]\right\} \right)\right]\leq\exp\left(\frac{\lambda^{2}\mathcal{L}^{2}}{2\varphi}\right).\label{eq:mgfg}
\end{equation}
As a consequence, }

\textit{
\begin{eqnarray*}
\mathbb{P}\left\{ f(U)-\mathbb{E}\left[f(U)\right]\leq-t\right\}  & \leq & \exp\left(-\frac{\varphi t^{2}}{2\mathcal{L}^{2}}\right),\\
\mathbb{P}\left\{ f(U)-\mathbb{E}\left[f(U)\right]\geq t\right\}  & \leq & \exp\left(-\frac{\varphi t^{2}}{2\mathcal{L}^{2}}\right).
\end{eqnarray*}
}\textbf{Remarks}. The proof involves the so-called ``inf-convolution''
argument and an application of the Brunn-Minkowski inequality; see
\cite{bl00} and \cite{m91}.\\
 \\
\textbf{\textit{Lemma A.2}}\textit{.}\textbf{\textit{ }}\textit{Assume
$U=\left\{ U_{i}\right\} _{i=1}^{n}$ consists of independent random
variables, all of which are supported on $[a,\,b]$. If $f:\,\mathbb{R}^{n}\rightarrow\mathbb{R}$
is separately convex}\footnote{Let the function $f_{j}:\,\mathbb{R}\rightarrow\mathbb{R}$ be defined
by varying only the $j$th co-ordinate of a function $f:\,\mathbb{R}^{n}\rightarrow\mathbb{R}$;
$f$ is \textit{separately convex} if for each $j\in\left\{ 1,\,2,\,...,\,n\right\} $,
$f_{j}$ is a convex function of the $j$th coordinate.}\textit{ and $\mathcal{L}-$Lipschitz with respect to the Euclidean
norm, then for all $\lambda\in\mathbb{R}$, 
\begin{equation}
\mathbb{E}\left[\exp\left(\lambda\left\{ f(U)-\mathbb{E}\left[f(U)\right]\right\} \right)\right]\leq\exp\left[\frac{\lambda^{2}(b-a)^{2}\mathcal{L}^{2}}{2}\right].\label{eq:mgfb}
\end{equation}
As a consequence, }

\textit{
\begin{eqnarray*}
\mathbb{P}\left[f(U)-\mathbb{E}\left[f(U)\right]\leq-t\right] & \leq & \exp\left(-\frac{t^{2}}{2\mathcal{L}^{2}(b-a)^{2}}\right),\\
\mathbb{P}\left[f(U)-\mathbb{E}\left[f(U)\right]\geq t\right] & \leq & \exp\left(-\frac{t^{2}}{2\mathcal{L}^{2}(b-a)^{2}}\right).
\end{eqnarray*}
}\textbf{Remarks}. One proof for Lemma A.2 involves the entropy method
and the so-called Herbst argument; see \cite{blm13}. Talagrand and
Ledoux have contributed to the result above in different papers.

\subsection{Proof of Propositions 2.1 and A.2.1 \label{subsec:Proof-of-Propositions 2.1 and A2.1}}

For any $q\in[1,\,\infty]$, $\left\Vert \frac{1}{n}X^{T}W\right\Vert _{q}$
is Lipschitz in $W$ with respect to the Euclidean norm. To see this,
note that a triangle inequality and a Cauchy-Schwarz inequality yield
\begin{eqnarray}
\left|\left\Vert \frac{1}{n}X^{T}W\right\Vert _{q}-\left\Vert \frac{1}{n}X^{T}W^{'}\right\Vert _{q}\right| & \leq & \left\Vert \frac{1}{n}X^{T}\left(W-W^{'}\right)\right\Vert _{q}\nonumber \\
 & \leq & \frac{1}{\sqrt{n}}\left\Vert \sqrt{\frac{1}{n}\sum_{i=1}^{n}X_{i}^{2}}\right\Vert _{q}\left\Vert W-W^{'}\right\Vert _{2}.\label{eq:lip}
\end{eqnarray}
As a result of Lemma A.1, we have the concentration in (\ref{eq:prop2-1}).

If $h(\theta^{*})\in\Omega$, (\ref{eq:prop2-1}) then implies that
$\left(\theta^{*},\,\mathbf{0}_{p}\right)$ ($\left(\theta^{*},\,0\right)$)
is an optimal solution to (\ref{eq:p1-1}) (respectively, (\ref{eq:p2-1})).
If $h(\theta^{*})\notin\Omega$, since $\left\{ \theta\in\mathbb{R}^{p}\,:\,h(\theta)\in\Omega\right\} \neq\emptyset$,
we can find some $\tilde{\theta}_{\alpha}$ such that $h(\tilde{\theta}_{\alpha})\in\Omega$.
Letting 
\begin{eqnarray*}
\tilde{\mu}_{\alpha} & = & \frac{1}{n}\sum_{i=1}^{n}X_{i}\left[Y_{i}-g\left(V_{i};\tilde{\theta}_{\alpha}\right)\right]-\frac{1}{n}\sum_{i=1}^{n}X_{i}\left[Y_{i}-g\left(V_{i};\theta^{*}\right)\right]\\
 & = & \frac{1}{n}\sum_{i=1}^{n}X_{i}\left[g\left(V_{i};\theta^{*}\right)-g\left(V_{i};\tilde{\theta}_{\alpha}\right)\right],
\end{eqnarray*}
(\ref{eq:prop2-1}) then implies that $\left(\tilde{\theta}_{\alpha},\,\tilde{\mu}_{\alpha}\right)$
is a feasible solution to (\ref{eq:p1-1}) with probability at least
$1-\alpha$. As a result, $\left(\tilde{\theta}_{\alpha},\,\left\Vert \tilde{\mu}_{\alpha}\right\Vert _{q}\right)$
is also a feasible solution to (\ref{eq:p2-1}) with probability at
least $1-\alpha$. 

In any case, an optimal solution $\left(\hat{\theta}_{\alpha}^{*},\,\hat{\mu}_{\alpha}^{*}\right)$
to (\ref{eq:p1-1}) must satisfy 
\begin{align*}
 & \left\Vert \frac{1}{n}\sum_{i=1}^{n}X_{i}\left[Y_{i}-g\left(V_{i};\hat{\theta}_{\alpha}^{*}\right)\right]-\frac{1}{n}\sum_{i=1}^{n}X_{i}\left[Y_{i}-g\left(V_{i};\theta^{*}\right)\right]\right\Vert _{\tilde{q}}\\
= & \left\Vert \frac{1}{n}\sum_{i=1}^{n}X_{i}\left[g\left(V_{i};\theta^{*}\right)-g\left(V_{i};\hat{\theta}_{\alpha}^{*}\right)\right]\right\Vert _{\tilde{q}}\geq\left\Vert \hat{\mu}_{\alpha}^{*}\right\Vert _{\tilde{q}}
\end{align*}
with probability at least $1-\alpha$.\textcolor{red}{{} }Similarly,
an optimal solution $\left(\hat{\theta}_{\alpha}^{*},\,\hat{\mu}_{\alpha}^{*}\right)$
to (\ref{eq:p2-1}) must satisfy 
\begin{align*}
 & \left\Vert \frac{1}{n}\sum_{i=1}^{n}X_{i}\left[Y_{i}-g\left(V_{i};\hat{\theta}_{\alpha}^{*}\right)\right]-\frac{1}{n}\sum_{i=1}^{n}X_{i}\left[Y_{i}-g\left(V_{i};\theta^{*}\right)\right]\right\Vert _{q}\\
= & \left\Vert \frac{1}{n}\sum_{i=1}^{n}X_{i}\left[g\left(V_{i};\theta^{*}\right)-g\left(V_{i};\hat{\theta}_{\alpha}^{*}\right)\right]\right\Vert _{q}\geq\hat{\mu}_{\alpha}^{*}
\end{align*}
with probability at least $1-\alpha$. 

On the other hand, in terms of (\ref{eq:p1-1}), applying the triangle
inequality yields 
\begin{align*}
 & \left\Vert \frac{1}{n}\sum_{i=1}^{n}X_{i}\left[g\left(V_{i};\theta^{*}\right)-g\left(V_{i};\hat{\theta}_{\alpha}^{*}\right)\right]-\hat{\mu}_{\alpha}^{*}\right\Vert _{q}\\
\leq & \left\Vert \frac{1}{n}X^{T}W\right\Vert _{q}+\left\Vert \frac{1}{n}\sum_{i=1}^{n}X_{i}\left[Y_{i}-g\left(V_{i};\hat{\theta}_{\alpha}^{*}\right)\right]-\hat{\mu}_{\alpha}^{*}\right\Vert _{q}\leq2r_{\alpha,q}^{*}
\end{align*}
with probability at least $1-\alpha$. In terms of (\ref{eq:p2-1}),
we simply have 
\[
\mathbb{P}\left(\left\Vert \frac{1}{n}\sum_{i=1}^{n}X_{i}\left[g\left(V_{i};\theta^{*}\right)-g\left(V_{i};\hat{\theta}_{\alpha}^{*}\right)\right]\right\Vert _{q}\leq2r_{\alpha,q}^{*}+\hat{\mu}_{\alpha}^{*}\right)\geq1-\alpha.
\]

\subsection{Proof of Theorem 2.1\label{subsec:Proof-of-Theorem 2.1}}

We have already derived (\ref{eq:13-2}) in Section \ref{sec:Gaussian-Regressions}.
To show (\ref{eq:13-1}), we define the event 
\[
\mathcal{E}=\left\{ \frac{\sigma}{R}\sum_{r=1}^{R}\left\Vert \frac{1}{n}X^{T}Z_{r}\right\Vert _{q}\geq\mathbb{E}_{W}\left[\left\Vert \frac{1}{n}X^{T}W\right\Vert _{q}\right]+\sqrt{\frac{1}{R}}\tau_{\beta_{1},q}\right\} .
\]
Like we have argued for (\ref{eq:c2}), we also have the upper deviation
inequality 
\begin{equation}
\mathbb{P}\left\{ \frac{\sigma}{R}\sum_{r=1}^{R}\left\Vert \frac{1}{n}X^{T}Z_{r}\right\Vert _{q}\geq\mathbb{E}_{W}\left[\left\Vert \frac{1}{n}X^{T}W\right\Vert _{q}\right]+t\right\} \leq\exp\left(\frac{-nRt^{2}}{2\sigma^{2}\left\Vert \sqrt{\frac{1}{n}\sum_{i=1}^{n}X_{i}^{2}}\right\Vert _{q}^{2}}\right)\label{eq:c2_u}
\end{equation}
and consequently, $\mathbb{P}\left(\mathcal{E}\right)\leq\beta_{1}$.
Let $\mathcal{E}^{c}$ denote the complement of $\mathcal{E}$. Under
$H_{1}$ and (\ref{eq:ms}), we have 
\begin{align*}
 & \mathbb{P}\left\{ \Psi_{q}(\hat{\theta}_{\alpha})\leq r_{\alpha,q}\right\} \\
= & \mathbb{P}\left\{ \Psi_{q}(\hat{\theta}_{\alpha})\leq r_{\alpha,q}\vert\mathcal{E}^{c}\right\} \mathbb{P}\left(\mathcal{E}^{c}\right)+\mathbb{P}\left\{ \Psi_{q}(\hat{\theta}_{\alpha})\leq r_{\alpha,q}\vert\mathcal{E}\right\} \mathbb{P}\left(\mathcal{E}\right)\\
\leq & \mathbb{P}\left\{ \Psi_{q}(\hat{\theta}_{\alpha})\leq r_{\alpha,q}\vert\mathcal{E}^{c}\right\} +\mathbb{P}\left(\mathcal{E}\right)\\
\leq & \mathbb{P}\left\{ \left\Vert \frac{1}{n}\sum_{i=1}^{n}X_{i}\left[g\left(V_{i};\theta^{*}\right)-g\left(V_{i};\hat{\theta}_{\alpha}\right)\right]\right\Vert _{q}-\left\Vert \frac{1}{n}X^{T}W\right\Vert _{q}\leq r_{\alpha,q}\vert\mathcal{E}^{c}\right\} +\beta_{1}\\
\leq & \mathbb{P}\left\{ \delta_{\alpha,\beta,q}-\left\Vert \frac{1}{n}X^{T}W\right\Vert _{q}\leq r_{\alpha,q}\vert\mathcal{E}^{c}\right\} +\beta_{1}\\
\leq & \mathbb{P}\left\{ \left\Vert \frac{1}{n}X^{T}W\right\Vert _{q}\geq\mathbb{E}_{W}\left[\left\Vert \frac{1}{n}X^{T}W\right\Vert _{q}\right]+\tau_{\beta_{2},q}\vert\mathcal{E}^{c}\right\} +\beta_{1}\\
\leq & \beta.
\end{align*}
In the above, the fifth line follows from (\ref{eq:ms}); the sixth
line follows from (\ref{eq:ms1}), the fact that we are conditioning
on $\mathcal{E}^{c}$, as well as (\ref{eq:prop2-1}); the last line
follows from the fact that $W$ is independent of $Z$.

\subsection{Additional derivations\label{subsec:Additional-Derivations}}

To show (\ref{eq:lead}), we define an i.i.d. sequence of Gaussian
random variables 
\[
\widetilde{W}_{k}\sim\mathcal{N}\left(0,\,\min_{j,l\in\left\{ 1,...,L\right\} }\frac{1}{2n^{2}}\sum_{i=1}^{n}\left(X_{ij}-X_{il}\right)^{2}\right)
\]
for $k=1,...,L$. Note that we have 
\[
\mathbb{E}_{W}\left[\left(\frac{1}{n}X_{j}^{T}W-\frac{1}{n}X_{l}^{T}W\right)^{2}\right]\geq\mathbb{E}_{\widetilde{W}}\left[\left(\widetilde{W}_{j}-\widetilde{W}_{l}\right)^{2}\right].
\]
By the Sudakov-Fernique Gaussian comparison result (see Corollary
3.14 in \cite{lt91}), we obtain 
\begin{eqnarray*}
\mathbb{E}_{W}\left[\left\Vert \frac{1}{n}X^{T}W\right\Vert _{\infty}\right] & \geq & \mathbb{E}_{W}\left[\max_{j\in\left\{ 1,...,L\right\} }\frac{1}{n}X_{j}^{T}W\right]\\
 & \geq & \frac{1}{2}\mathbb{E}_{\widetilde{W}}\left[\max_{j\in\left\{ 1,...,L\right\} }\widetilde{W}_{j}\right]\\
 & \geq & \frac{1}{2}\left(1-\frac{1}{e}\right)\sqrt{\frac{\log L}{4n^{2}}\min_{j,l\in\left\{ 1,...,L\right\} }\sum_{i=1}^{n}\left(X_{ij}-X_{il}\right)^{2}}
\end{eqnarray*}
(for all $L\geq20$), where the last line follows from a classical
lower bound on the Gaussian maximum (see, e.g., \cite{lt91}). The
upper bound 
\[
\mathbb{E}_{W}\left[\left\Vert \frac{1}{n}X^{T}W\right\Vert _{\infty}\right]\leq\sqrt{\frac{2\log L}{n^{2}}\max_{j\in\left\{ 1,...,L\right\} }\sum_{i=1}^{n}X_{ij}^{2}}+\sqrt{\frac{8}{n^{2}\log L}\max_{j\in\left\{ 1,...,L\right\} }\sum_{i=1}^{n}X_{ij}^{2}}
\]
(for all $L\geq2$) is another existing result on the Gaussian maximum
(see, e.g., \cite{w15}). \\
 \\
\textbf{Remarks}. To obtain the lower bound on $\mathbb{E}_{W}\left[\left\Vert \frac{1}{n}X^{T}W\right\Vert _{\infty}\right]$,
we first compare the dependent sequence $\left\{ \frac{1}{n}X_{j}^{T}W\right\} _{j=1}^{L}$
with another independent Gaussian sequence $\widetilde{W}=\left\{ \widetilde{W}_{j}\right\} _{j=1}^{L}$
and then apply a lower bound on $\mathbb{E}_{\widetilde{W}}\left[\max_{j\in\left\{ 1,...,L\right\} }\widetilde{W}_{j}\right]$.
In contrast, the upper bound on $\mathbb{E}_{W}\left[\left\Vert \frac{1}{n}X^{T}W\right\Vert _{\infty}\right]$
is obtained by applying $\sum_{j=1}^{L}\mathbb{P}\left(\left|\frac{1}{n}X_{j}^{T}W\right|\geq t\right)$,
where independence is not needed. Moreover, the upper bound also holds
when $W$ is a sequence of sub-Gaussian variables while the lower
bound requires $W$ to be a sequence of Gaussian variables.

\subsection{Proofs of Lemmas 4.1 and 4.2}

As a result of Lemma A.1 and (\ref{eq:lip}), we have the concentration
in Lemma 4.1. Because $\left\Vert \frac{1}{n}X^{T}W\right\Vert _{q}$
is separately convex in terms of $W$, Lemma A.2 implies the concentration
in Lemma 4.2.

\subsection{Proof of Proposition 4.1 }

Using the argument that leads to (\ref{eq:lip}), we can show $\left\Vert \frac{1}{n}\sum_{i=1}^{n}X_{i}\left[Y_{i}-g\left(V_{i};\theta^{*}\right)\right]\right\Vert _{q}$
is Lipschitz in $Y$ with respect to the Euclidean norm for any $q\in[1,\,\infty]$.
That is, 
\begin{align}
 & \left|\left\Vert \frac{1}{n}\sum_{i=1}^{n}X_{i}\left[Y_{i}-g\left(V_{i};\theta^{*}\right)\right]\right\Vert _{q}-\left\Vert \frac{1}{n}\sum_{i=1}^{n}X_{i}\left[Y_{i}^{'}-g\left(V_{i};\theta^{*}\right)\right]\right\Vert _{q}\right|\nonumber \\
\leq & \frac{1}{\sqrt{n}}\left\Vert \sqrt{\frac{1}{n}\sum_{i=1}^{n}X_{i}^{2}}\right\Vert _{q}\left\Vert Y-Y^{'}\right\Vert _{2}.\label{eq:lip1}
\end{align}
Note that $\left\Vert \frac{1}{n}\sum_{i=1}^{n}X_{i}\left[Y_{i}-g\left(V_{i};\theta^{*}\right)\right]\right\Vert _{q}$
is separately convex in terms of $Y$. As a result of Lemma A.2, we
have the concentration in (\ref{eq:prop3.1}).

To establish (\ref{eq:sandwich}), we exploit the convexity of $l_{q}-$norms
and the fact that $\mathbb{E}\left(Y_{i}|V_{i}\right)=g\left(V_{i};\theta^{*}\right)$.
Let $Y^{'}=\left\{ Y_{i}^{'}\right\} _{i=1}^{n}$ be an independent
sequence identical to but independent of $Y$ conditioning on $V$,
and $\varepsilon=\left\{ \varepsilon_{i}\right\} _{i=1}^{n}$ be i.i.d.
Radamacher random variables independent of $Y$, $Y^{'}$, and $V$.
We obtain 
\begin{align}
 & \mathbb{E}_{Y|V}\left\{ \left\Vert \frac{1}{n}\sum_{i=1}^{n}X_{i}\left[Y_{i}-g\left(V_{i};\theta^{*}\right)\right]\right\Vert _{q}\right\} \nonumber \\
= & \mathbb{E}_{Y|V}\left\{ \left\Vert \frac{1}{n}\sum_{i=1}^{n}X_{i}\left[Y_{i}-\mathbb{E}_{Y_{i}^{'}\vert V_{i}}\left(Y_{i}^{'}\right)\right]\right\Vert _{q}\right\} \nonumber \\
= & \mathbb{E}_{Y|V}\left\{ \left\Vert \mathbb{E}_{Y^{'}\vert V}\left[\frac{1}{n}\sum_{i=1}^{n}X_{i}\left(Y_{i}-Y_{i}^{'}\right)\right]\right\Vert _{q}\right\} \nonumber \\
\leq & \mathbb{E}_{Y^{'},Y|V}\left\{ \left\Vert \frac{1}{n}\sum_{i=1}^{n}X_{i}\left(Y_{i}-Y_{i}^{'}\right)\right\Vert _{q}\right\} \nonumber \\
= & \mathbb{E}_{\varepsilon,Y^{'},Y|V}\left\{ \left\Vert \frac{1}{n}\sum_{i=1}^{n}\varepsilon_{i}X_{i}\left(Y_{i}-Y_{i}^{'}\right)\right\Vert _{q}\right\} \nonumber \\
\leq & 2\mathbb{E}_{\varepsilon,Y|V}\left\{ \left\Vert \frac{1}{n}\sum_{i=1}^{n}\varepsilon_{i}X_{i}Y_{i}\right\Vert _{q}\right\} ,\label{eq:sym}
\end{align}
where the second line follows since $\mathbb{E}\left(Y_{i}^{'}|V_{i}\right)=g\left(V_{i};\theta^{*}\right)$,
the fourth line follows from Jensen's inequality, and the sixth line
follows from the fact that $\varepsilon_{i}X_{i}\left(Y_{i}-Y_{i}^{'}\right)$
and $X_{i}\left(Y_{i}-Y_{i}^{'}\right)$ have the same distribution.

On the other hand, similar argument from above also yields 
\begin{align*}
 & \mathbb{E}_{\varepsilon,Y|V}\left\{ \left\Vert \frac{1}{2n}\sum_{i=1}^{n}\varepsilon_{i}X_{i}\left[Y_{i}-g\left(V_{i};\theta^{*}\right)\right]\right\Vert _{q}\right\} \\
= & \mathbb{E}_{\varepsilon,Y|V}\left\{ \left\Vert \frac{1}{2n}\sum_{i=1}^{n}\varepsilon_{i}X_{i}\left[Y_{i}-\mathbb{E}_{Y_{i}^{'}\vert V_{i}}\left(Y_{i}^{'}\right)\right]\right\Vert _{q}\right\} \\
\leq & \mathbb{E}_{\varepsilon,Y^{'},Y|V}\left\{ \left\Vert \frac{1}{2n}\sum_{i=1}^{n}\varepsilon_{i}X_{i}\left(Y_{i}-Y_{i}^{'}\right)\right\Vert _{q}\right\} \\
= & \mathbb{E}_{Y^{'},Y|V}\left\{ \left\Vert \frac{1}{2n}\sum_{i=1}^{n}X_{i}\left(Y_{i}-Y_{i}^{'}\right)\right\Vert _{q}\right\} .
\end{align*}
Applying the following inequality 
\begin{align*}
 & \left\Vert \frac{1}{2n}\sum_{i=1}^{n}X_{i}\left(Y_{i}-Y_{i}^{'}\right)\right\Vert _{q}\\
\leq & \left\Vert \frac{1}{2n}\sum_{i=1}^{n}X_{i}\left(Y_{i}-g\left(V_{i};\theta^{*}\right)\right)\right\Vert _{q}+\left\Vert \frac{1}{2n}\sum_{i=1}^{n}X_{i}\left(Y_{i}^{'}-g\left(V_{i};\theta^{*}\right)\right)\right\Vert _{q},
\end{align*}
and taking expectations gives 
\[
\mathbb{E}_{Y^{'},Y|V}\left\{ \left\Vert \frac{1}{2n}\sum_{i=1}^{n}X_{i}\left(Y_{i}-Y_{i}^{'}\right)\right\Vert _{q}\right\} \leq\mathbb{E}_{Y\vert V}\left\{ \left\Vert \frac{1}{n}\sum_{i=1}^{n}X_{i}\left(Y_{i}-g\left(V_{i};\theta^{*}\right)\right)\right\Vert _{q}\right\} .
\]
Putting the pieces together, we obtain the result in (\ref{eq:sandwich}).

\subsection{Proof of Proposition 4.2}

We first show that $\mathbb{E}_{\varepsilon}\left\{ \left\Vert \frac{1}{n}\sum_{i=1}^{n}\varepsilon_{i}X_{i}Y_{i}\right\Vert _{q}\right\} $
is Lipschitz in $Y$ with respect to the Euclidean norm for any $q\in[1,\,\infty]$.
That is, 
\begin{align*}
 & \left|\mathbb{E}_{\varepsilon}\left\{ \left\Vert \frac{1}{n}\sum_{i=1}^{n}\varepsilon_{i}X_{i}Y_{i}\right\Vert _{q}\right\} -\mathbb{E}_{\varepsilon}\left\{ \left\Vert \frac{1}{n}\sum_{i=1}^{n}\varepsilon_{i}X_{i}Y_{i}^{'}\right\Vert _{q}\right\} \right|\\
\leq & \frac{1}{\sqrt{n}}\left\Vert \sqrt{\frac{1}{n}\sum_{i=1}^{n}X_{i}^{2}}\right\Vert _{q}\sqrt{\mathbb{E}_{\varepsilon}\left[\sum_{i=1}^{n}\varepsilon_{i}^{2}\left(Y_{i}-Y_{i}^{'}\right)^{2}\right]}\\
\leq & \frac{1}{\sqrt{n}}\left\Vert \sqrt{\frac{1}{n}\sum_{i=1}^{n}X_{i}^{2}}\right\Vert _{q}\left\Vert Y-Y^{'}\right\Vert _{2}.
\end{align*}
Note that $\mathbb{E}_{\varepsilon}\left\{ \left\Vert \frac{1}{n}\sum_{i=1}^{n}\varepsilon_{i}X_{i}Y_{i}\right\Vert _{q}\right\} $
is separately convex in terms of $Y$. As a result of Lemma A.2, we
have the following concentration 
\begin{equation}
\mathbb{P}\left\{ \mathbb{E}_{\varepsilon,Y\vert V}\left\{ \left\Vert \frac{1}{n}\sum_{i=1}^{n}\varepsilon_{i}X_{i}Y_{i}\right\Vert _{q}\right\} \geq\mathbb{E}_{\varepsilon}\left\{ \left\Vert \frac{1}{n}\sum_{i=1}^{n}\varepsilon_{i}X_{i}Y_{i}\right\Vert _{q}\right\} +\tau_{\alpha_{2},q}\right\} \leq\alpha_{2}.\label{eq:c1-3}
\end{equation}

Let $\varepsilon=\left\{ \varepsilon_{i}\right\} _{i=1}^{n}$ be an
i.i.d. sequence of Radamacher random variables, independent of $Y$
and $V$. We can again show that $\left\Vert \frac{1}{n}\sum_{i=1}^{n}\varepsilon_{i}Y_{i}X_{i}\right\Vert _{q}$
is Lipschitz in $\varepsilon$ with respect to the Euclidean norm
for any $q\in[1,\,\infty]$ and the Lipschitz constant\footnote{Like $\left\Vert \sqrt{\frac{1}{n}\sum_{i=1}^{n}X_{i}^{2}}\right\Vert _{q}$,
we define 
\begin{eqnarray*}
\left\Vert \sqrt{\frac{1}{n}\sum_{i=1}^{n}\left(Y_{i}X_{i}\right)^{2}}\right\Vert _{q} & = & \sqrt[q]{\sum_{j=1}^{L}\left(\sqrt{\frac{1}{n}\sum_{i=1}^{n}Y_{i}^{2}X_{ij}^{2}}\right)^{q}},\qquad q\in[1,\,\infty)\\
\left\Vert \sqrt{\frac{1}{n}\sum_{i=1}^{n}\left(Y_{i}X_{i}\right)^{2}}\right\Vert _{q} & = & \max_{j\in\left\{ 1,...,L\right\} }\sqrt{\frac{1}{n}\sum_{i=1}^{n}Y_{i}^{2}X_{ij}^{2}},\qquad q=\infty.
\end{eqnarray*}
} is $\frac{1}{\sqrt{n}}\left\Vert \sqrt{\frac{1}{n}\sum_{i=1}^{n}\left(Y_{i}X_{i}\right)^{2}}\right\Vert _{q}$,
which is bounded from above by $\frac{1}{\sqrt{n}}\left\Vert \sqrt{\frac{1}{n}\sum_{i=1}^{n}X_{i}^{2}}\right\Vert _{q}$
given $0\leq Y_{i}\leq1$. Let $\left\{ \varepsilon_{ir}:\,i=1,...,n,\,r=1,...,R\right\} $
be a collection of i.i.d. Radamacher random draws, independent of
$Y$ and $V$. Conditioning on $Y$ and $V$, (\ref{eq:subn}) and
(\ref{eq:mgfb}) imply $\frac{1}{R}\sum_{r=1}^{R}\left\Vert \frac{1}{n}\sum_{i=1}^{n}\varepsilon_{ir}Y_{i}X_{i}\right\Vert _{q}-\mathbb{E}_{\varepsilon}\left(\left\Vert \frac{1}{n}\sum_{i=1}^{n}\varepsilon_{i}X_{i}Y_{i}\right\Vert _{q}\right)$
is sub-Gaussian with parameter at most $\frac{2}{\sqrt{nR}}\left\Vert \sqrt{\frac{1}{n}\sum_{i=1}^{n}X_{i}^{2}}\right\Vert _{q}$.
Therefore, we have 
\begin{align*}
 & \mathbb{E}_{Y\vert V}\left[\mathbb{E}_{\varepsilon}\left[\exp\left(\lambda\left[\frac{1}{R}\sum_{r=1}^{R}\left\Vert \frac{1}{n}\sum_{i=1}^{n}\varepsilon_{ir}Y_{i}X_{i}\right\Vert _{q}-\mathbb{E}_{\varepsilon}\left(\left\Vert \frac{1}{n}\sum_{i=1}^{n}\varepsilon_{i}X_{i}Y_{i}\right\Vert _{q}\right)\right]\right)\right]\right]\\
\leq & \exp\left[\lambda^{2}\frac{4\left\Vert \sqrt{\frac{1}{n}\sum_{i=1}^{n}X_{i}^{2}}\right\Vert _{q}^{2}}{2nR}\right].
\end{align*}
Consequently, (\ref{eq:lower}) yields the following concentration
\begin{equation}
\mathbb{P}\left\{ \mathbb{E}_{\varepsilon}\left\{ \left\Vert \frac{1}{n}\sum_{i=1}^{n}\varepsilon_{i}X_{i}Y_{i}\right\Vert _{q}\right\} \geq\frac{1}{R}\sum_{r=1}^{R}\left\Vert \frac{1}{n}\sum_{i=1}^{n}\varepsilon_{ir}Y_{i}X_{i}\right\Vert _{q}+\frac{2}{\sqrt{R}}\tau_{\alpha_{3},q}\right\} \leq\alpha_{3}\label{eq:c2-1}
\end{equation}
Combining (\ref{eq:prop3.1}), (\ref{eq:sym}), (\ref{eq:c1-3}) and
(\ref{eq:c2-1}) yields (\ref{eq:c1-2}).

\end{document}